\documentclass[11pt]{article}
\usepackage{amscd}
\usepackage{amssymb,amsfonts,amsmath,psfrag,amscd,cite}
\usepackage{amsmath}
  \usepackage{paralist}
  \usepackage{graphics}
  \usepackage{epsfig}
  \usepackage[colorlinks=true]{hyperref}
  \usepackage[all]{xy}
\usepackage{graphicx}
\newtheorem{theo}{Theorem}

\newtheorem{rema}[theo]{Remark}

\makeatletter \@addtoreset{equation}{section}
\@addtoreset{theo}{section} \makeatother

\begin{document}
\date{}
\title{Symmetric Reduction and Hamilton-Jacobi Equation of\\ Rigid Spacecraft with a Rotor}
\author{Hong Wang  \\
School of Mathematical Sciences and LPMC,\\
Nankai University, Tianjin 300071, P.R.China\\
E-mail: hongwang@nankai.edu.cn\\
January 1, 2014} \maketitle

{\bf Abstract.} In this paper, we consider the rigid spacecraft with
an internal rotor as a regular point reducible regular controlled
Hamiltonian (RCH) system. In the cases of coincident and
non-coincident centers of buoyancy and gravity, we give explicitly
the motion equation and Hamilton-Jacobi equation of reduced
spacecraft-rotor system on a symplectic leaf by calculation in
detail, respectively, which show the effect on controls
in regular symplectic reduction and Hamilton-Jacobi theory.\\



\tableofcontents

\section{Introduction}

It is well-known that the theory of controlled mechanical systems
has formed an important subject in recent years. Its research
gathers together some separate areas of research such as mechanics,
differential geometry and nonlinear control theory, etc., and the
emphasis of this research on geometry is motivated by the aim of
understanding the structure of equations of motion of the system in
a way that helps both analysis and design. There is a natural
two-fold method in the study of controlled mechanical systems.
First, it is a special class of nonlinear control systems
considered, and one can study it by using the feedback control and
optimal control method. The second, it is also a special class of
mechanical systems, one can study it combining with the analysis of
dynamic systems and the geometric reduction theory of Hamiltonian
and Lagrangian systems. Thus, the theory of controlled mechanical
systems presents a challenging and promising research area between
the study of classical mechanics and that of modern nonlinear
geometric control theory, and a lot of researchers are absorbed to
pour into the area and there have been a lot of interesting results.
Such as, Bloch et al in \cite{blkrmaal92, blle02, bllema00,
bllema01}, studied the symmetry and control as well as referred to
the use of feedback control to realize a modification to the
structure of a given mechanical system; Nijmeijer and Van der Schaft
in \cite{nivds90}, studied the nonlinear dynamical control systems
as well as the use of feedback control to stabilize mechanical
systems, and Van der Schaft in \cite{vds00,vds06} referred to the
reduction and control of implicit (port) Hamiltonian systems.\\

In particular, we note that in Marsden et al \cite{mawazh10}, the
authors studied regular reduction theory of controlled Hamiltonian
systems with symplectic structure and symmetry, as an extension of
regular symplectic reduction theory of Hamiltonian systems under
regular controlled Hamiltonian equivalence conditions. Wang in
\cite{wa12} generalized the work in \cite{mawazh10} to study the
singular reduction theory of regular controlled Hamiltonian systems,
and Wang and Zhang in \cite{wazh12} generalized the work in
\cite{mawazh10} to study optimal reduction theory of controlled
Hamiltonian systems with Poisson structure and symmetry by using
optimal momentum map and reduced Poisson tensor (or reduced
symplectic form), and Ratiu and Wang in \cite{rawa12} studied the
Poisson reduction of controlled Hamiltonian system by
controllability distribution. These research work not only gave a
variety of reduction methods for controlled Hamiltonian systems, but
also showed a variety of relationships of controlled Hamiltonian
equivalence of these systems.\\

At the same time, we note also that Hamilton-Jacobi theory is an
important part of classical mechanics. On the one hand, it provides
a characterization of the generating functions of certain
time-dependent canonical transformations, such that a given
Hamiltonian system in such a form that its solutions are extremely
easy to find by reduction to the equilibrium, see Abraham and
Marsden \cite{abma78}, Arnold \cite{ar89} and Marsden and Ratiu
\cite{mara99}. On the other hand, it is possible in many cases that
Hamilton-Jacobi theory provides an immediate way to integrate the
equation of motion of system, even when the problem of Hamiltonian
system itself has not been or cannot be solved completely. In
addition, the Hamilton-Jacobi equation is also fundamental in the
study of the quantum-classical relationship in quantization, and it
also plays an important role in the development of numerical
integrators that preserve the symplectic structure and in the study
of stochastic dynamical systems, see Woodhouse \cite{wo92}, Ge and
Marsden \cite{gema88}, Marsden and West \cite{mawe01} and
L\'{a}zaro-Cam\'{i} and Ortega \cite{laor09}. For these reasons it
is described as a useful tools in the study of Hamiltonian system
theory, and has been extensively developed in past many years and
become one of the most active subjects in the study of modern
applied mathematics and analytical mechanics, which absorbed a lot
of researchers to pour into it and a lot of deep and beautiful
results have been obtained, see Cari$\tilde{n}$ena et al
\cite{cagrmamamuro06} and \cite{cagrmamamuro10}, Iglesias et
al\cite{iglema08}, Le\'{o}n et al \cite{balemamamu12, lemama09,
lemama10} Ohsawa and Bloch \cite{ohbl09}
for more details.\\

Just as we have known that Hamilton-Jacobi theory from the
variational point of view is originally developed by Jacobi in 1866,
which state that the integral of Lagrangian of a system along the
solution of its Euler-Lagrange equation satisfies the
Hamilton-Jacobi equation. The classical description of this problem
from the geometrical point of view is given by Abraham and Marsden
in \cite{abma78}, and it is developed in the context of
time-dependent Hamiltonian system by Marsden and Ratiu in
\cite{mara99}. The Hamilton-Jacobi equation may be regarded as a
nonlinear partial differential equation for some generating function
$S$, and the problem is become how to choose a time-dependent
canonical transformation $\Psi: T^*Q\times \mathbb{R} \rightarrow
T^*Q\times \mathbb{R}, $ which transforms the dynamical vector field
of a time-dependent Hamiltonian system to equilibrium, such that the
generating function $S$ of $\Psi$ satisfies the time-dependent
Hamilton-Jacobi equation. In particular, for the time-independent
Hamiltonian system, we may look for a symplectic map as the
canonical transformation. This work offers a important idea that one
can use the dynamical vector field of a Hamiltonian system to
describe Hamilton-Jacobi equation. Moreover, assume that $\gamma: Q
\rightarrow T^*Q$ is a closed one-form on smooth manifold $Q$, and
define that $X_H^\gamma = T\pi_{Q}\cdot X_H \cdot \gamma$, where
$X_{H}$ is the dynamical vector field of Hamiltonian system
$(T^*Q,\omega,H)$. Then the fact that $X_H^\gamma$ and $X_H$ are
$\gamma$-related, that is, $T\gamma\cdot X_H^\gamma= X_H\cdot
\gamma$ is equivalent that $\mathbf{d}(H \cdot \gamma)=0, $ which is
given in Cari\~{n}ena et al \cite{cagrmamamuro06} and Iglesias et
al\cite{iglema08}. Since the Hamilton-Jacobi theory is developed
based on the Hamiltonian picture of dynamics, motivated by the above
research work, Wang in \cite{wa13a} uses the dynamical vector field
of Hamiltonian system and the regular reduced Hamiltonian system to
describe the Hamilton-Jacobi theory for these systems. Moreover,
Wang in \cite{wa13d} extends the Hamilton-Jacobi theory to the
regular controlled Hamiltonian system and its regular reduced
systems, and describes the relationship between the RCH-equivalence
for RCH systems and the solutions of corresponding Hamilton-Jacobi
equations.\\

Now, it is a natural problem if there is a practical controlled
Hamiltonian system and how to show the effect on controls in regular
symplectic reduction and Hamilton-Jacobi theory of the system. In
this paper, as an application of the regular point symplectic
reduction and Hamilton-Jacobi theory of RCH system with symmetry, we
consider that rigid spacecraft with an internal rotor is a regular
point reducible RCH system given in Marsden et al \cite{mawazh10},
where the rigid spacecraft with an internal rotor is modeled as a
Hamiltonian system with control, which is given in Bloch and Leonard
\cite{blle02}, Bloch, Leonard and Marsden \cite{bllema01}. In the
cases of coincident and non-coincident centers of buoyancy and
gravity, we give explicitly the motion equation and Hamilton-Jacobi
equation of reduced spacecraft-rotor system on a symplectic leaf by
calculation in detail, respectively. These equations are more
complex than that of Hamiltonian system without control and describe
explicitly the effect on controls
in regular symplectic reduction and Hamilton-Jacobi theory.\\

A brief of outline of this paper is as follows. In the second
section, we first review some relevant definitions and basic facts
about rigid spacecraft with an internal rotor, which will be used in
subsequent sections. As an application of the theoretical result of
symplectic reduction of RCH system given by Marsden et al
\cite{mawazh10}, in the third section we regard the rigid spacecraft
with an internal rotor as a regular point reducible RCH system on
the generalization of rotation group $\textmd{SO}(3)\times S^1$ and
on the generalization of Euclidean group $\textmd{SE}(3)\times S^1$,
respectively, and in the cases of coincident and non-coincident
centers of buoyancy and gravity, we give explicitly the motion
equations of their reduced RCH systems on the symplectic leaves by
calculation in detail. Moreover, as an application of the
theoretical result of Hamilton-Jacobi theory of regular reduced RCH
system given by Wang \cite{wa13d}, in the fourth section, we give
the Hamilton-Jacobi equations of the reduced rigid spacecraft-rotor
systems on the symplectic leaves by calculation in detail,
respectively, in the cases of coincident and non-coincident centers
of buoyancy and gravity. These research work develop the application
of symplectic reduction and Hamilton-Jacobi theory of RCH systems
with symmetry and make us have much deeper understanding and
recognition for the structure of Hamiltonian systems and RCH
systems.

\section{The Rigid Spacecraft with a Rotor}

In this paper, our goal is to give the regular point reduction and
Hamilton-Jacobi theorem of rigid spacecraft with an internal rotor.
In order to do these, in this section, we first review some relevant
definitions and basic facts about rigid spacecraft with an internal
rotor, which will be used in subsequent sections. We shall follow
the notations and conventions introduced in Bloch and Leonard
\cite{blle02}, Bloch, Leonard and Marsden \cite{bllema01}, Marsden
\cite{ma92}, Marsden and Ratiu \cite{mara99}, and Marsden et al
\cite{mawazh10}. In this paper, we assume that all manifolds are
real, smooth and finite dimensional and all actions are smooth left
actions. For convenience, we also assume that all controls appearing
in this paper are the admissible controls.

\subsection{The Spacecraft-Rotor System with Coincident Centers}

We consider a rigid spacecraft (to be called the carrier body)
carrying an internal rotor, and assume that the only external forces
and torques acting on the spacecraft-rotor system are due to
buoyancy and gravity. In general, it is possible that the
spacecraft's center of buoyancy may not be coincident with its
center of gravity. But, in this subsection we assume that the
spacecraft is symmetric and to have uniformly distributed mass, and
the center of buoyancy and the center of gravity are coincident.
Denote by $O$ the center of mass of the system in the body frame and
at $O$ place a set of (orthogonal) body axes. Assume that the body
coordinate axes are aligned with principal axes of the carrier body,
and the rotor is aligned along the third principal axis, see Bloch
and Leonard \cite{blle02} and Bloch, Leonard and Marsden
\cite{bllema01}. The rotor spins under the influence of a torque $u$
acting on the rotor. If translations are ignored and only rotations
are considered, in this case, then the configuration space is
$Q=\textmd{SO}(3)\times S^1$, with the first factor being the
attitude of rigid spacecraft and the second factor being the angle
of rotor. The corresponding phase space is the cotangent bundle
$T^\ast Q =T^\ast \textmd{SO}(3)\times T^\ast S^1$, where $T^\ast
S^1 \cong T^\ast \mathbb{R}\cong \mathbb{R}\times \mathbb{R}$, with
the canonical symplectic form.\\

Let $I=diag(I_1,I_2,I_3)$ be the moment of inertia of the carrier
body in the principal body-fixed frame, and $J_3$ be the moment of
inertia of rotor around its rotation axis. Let $J_{3k},\; k=1,2,$ be
the moments of inertia of the rotor around the $k$th principal axis
with $k=1,2,$ and denote by $\bar{I}_k= I_k +J_{3k}, \; k=1,2, \;
\bar{I}_3= I_3. $ Let $\Omega=(\Omega_1,\Omega_2,\Omega_3)$ be the
vector of body angular velocities computed with respect to the axes
fixed in the body and $(\Omega_1,\Omega_2,\Omega_3)\in
\mathfrak{so}(3)$. Let $\alpha$ be the relative angle of rotor and
$\dot{\alpha}$ the rotor relative angular velocity about the third
principal axis with respect to a carrier body fixed frame. Consider
the Lagrangian of the system $L(A,\Omega,\alpha,\dot{\alpha}):
TQ\cong \textmd{SO}(3)\times\mathfrak{so}(3)\times
\mathbb{R}\times\mathbb{R}\to \mathbb{R}$, which is the total
kinetic energy of the rigid spacecraft plus the kinetic energy of
rotor, given by
$$L(A,\Omega,\alpha,\dot{\alpha})=\dfrac{1}{2}[\bar{I}_1\Omega_1^2
+\bar{I}_2\Omega_2^2+\bar{I}_3\Omega_3^2
+J_3(\Omega_3+\dot{\alpha})^2],$$ where $A\in \textmd{SO}(3)$,
$\Omega=(\Omega_1,\Omega_2,\Omega_3)\in \mathfrak{so}(3)$, $\alpha
\in \mathbb{R}$, $\dot{\alpha} \in \mathbb{R}$. If we introduce the
conjugate angular momentum, given by $\Pi_k= \dfrac{\partial
L}{\partial \Omega_k} =\bar{I}_k\Omega_k, \; k=1,2$, $\Pi_3=
\dfrac{\partial L}{\partial \Omega_3}
=\bar{I}_3\Omega_3+J_3(\Omega_3+\dot{\alpha}),\quad l =
\dfrac{\partial L}{\partial \dot{\alpha}}
=J_3(\Omega_3+\dot{\alpha}),$ and by the Legendre transformation
$FL:
\textmd{SO}(3)\times\mathfrak{so}(3)\times\mathbb{R}\times\mathbb{R}\to
\textmd{SO}(3)\times
\mathfrak{so}^\ast(3)\times\mathbb{R}\times\mathbb{R},\quad
(A,\Omega,\alpha,\dot{\alpha})\to(A,\Pi,\alpha,l)$, where
$\Pi=(\Pi_1,\Pi_2,\Pi_3)\in \mathfrak{so}^\ast(3)$, $l \in
\mathbb{R}$, we have the Hamiltonian $H(A,\Pi,\alpha,l): T^* Q \cong
\textmd{SO}(3)\times \mathfrak{so}^\ast
(3)\times\mathbb{R}\times\mathbb{R}\to \mathbb{R}$ given by
\begin{align}
H(A,\Pi,\alpha,l) &= \Omega\cdot \Pi+\dot{\alpha}\cdot
l-L(A,\Omega,\alpha,\dot{\alpha}) \nonumber \\
&= \frac{1}{2}[\frac{\Pi_1^2}{\bar{I}_1}+\frac{\Pi_2^2}{\bar{I}_2}
+\frac{(\Pi_3-l)^2}{\bar{I}_3}+\frac{l^2}{J_3}].
\end{align}

In order to give the motion equation of spacecraft-rotor system, in
the following we need to consider the symmetry and reduced
symplectic structure of the configuration space $Q=
\textmd{SO}(3)\times S^1 $.

\subsection{The Spacecraft-Rotor System with Non-coincident Centers}

Since it is possible that the spacecraft's center of buoyancy may
not be coincident with its center of gravity, in this subsection
then we consider the spacecraft-rotor system with non-coincident
centers of buoyancy and gravity. We fix an orthogonal coordinate
frame to the carrier body with origin located at the center of
buoyancy and axes aligned with the principal axes of the carrier
body, and the rotor is aligned along the third principal axis, see
Bloch and Leonard \cite{blle02}, Bloch, Leonard and Marsden
\cite{bllema01}, and Leonard and Marsden \cite{lema97}. The rotor
spins under the influence of a torque $u$ acting on the rotor. When
the carrier body is oriented so that the body-fixed frame is aligned
with the inertial frame, the third principal axis aligns with the
direction of gravity. The vector from the center of buoyancy to the
center of gravity with respect to the body-fixed frame is $h\chi$,
where $\chi$ is an unit vector on the line connecting the two
centers which is assumed to be aligned along the third principal
axis, and $h$ is the length of this segment. The mass of the carrier
body is denoted $m$ and the magnitude of gravitational acceleration
is $g$, and let $\Gamma$ be the unit vector viewed by an observer
moving with the body. In this case, the configuration space is
$Q=\textmd{SO}(3)\circledS \mathbb{R}^3\times S^1 \cong
\textmd{SE}(3)\times S^1$, with the first factor being the attitude
of rigid spacecraft and the drift of spacecraft in the rotational
process and the second factor being the angle of rotor. The
corresponding phase space is the cotangent bundle $T^\ast Q =T^\ast
\textmd{SE}(3)\times T^\ast S^1$, where $T^\ast S^1 \cong T^\ast
\mathbb{R}\cong \mathbb{R}\times \mathbb{R}$, with the
canonical symplectic form.\\

Consider the Lagrangian of the system
$L(A,c,\Omega,\Gamma,\alpha,\dot{\alpha}): TQ \cong
\textmd{SE}(3)\times\mathfrak{se}(3)\times\mathbb{R}\times\mathbb{R}\to
\mathbb{R}$, which is the total kinetic energy of the rigid
spacecraft plus the kinetic energy of rotor minus potential energy
of the system, given by
$$L(A,c,\Omega,\Gamma,\alpha,\dot{\alpha})=\dfrac{1}{2}[\bar{I}_1\Omega_1^2
+\bar{I}_2\Omega_2^2+\bar{I}_3\Omega_3^2
+J_3(\Omega_3+\dot{\alpha})^2]- mgh\Gamma \cdot \chi,$$ where
$(A,c)\in \textmd{SE}(3)$, $\Omega=(\Omega_1,\Omega_2,\Omega_3)\in
\mathfrak{so}(3)$, $\Gamma \in \mathbb{R}^3$, $(\Omega, \Gamma)\in
\mathfrak{se}(3)$, $\alpha \in \mathbb{R}$, $\dot{\alpha} \in
\mathbb{R}$. If we introduce the conjugate angular momentum, given
by $\Pi_k= \dfrac{\partial L}{\partial \Omega_k} =\bar{I}_k\Omega_k,
\; k=1,2 $, $\Pi_3= \dfrac{\partial L}{\partial \Omega_3}
=\bar{I}_3\Omega_3+J_3(\Omega_3+\dot{\alpha}),\quad l =
\dfrac{\partial L}{\partial \dot{\alpha}}
=J_3(\Omega_3+\dot{\alpha}),$ and by the Legendre transformation
$$FL:
\textmd{SE}(3)\times\mathfrak{se}(3)\times\mathbb{R}\times\mathbb{R}\to
\textmd{SE}(3)\times
\mathfrak{se}^\ast(3)\times\mathbb{R}\times\mathbb{R},\quad
(A,c,\Omega,\Gamma,\alpha,\dot{\alpha})\to(A,c,\Pi,\Gamma,\alpha,l),
$$ where $\Pi=(\Pi_1,\Pi_2,\Pi_3)\in \mathfrak{so}^\ast(3)$, $(\Pi,
\Gamma)\in \mathfrak{se}^\ast(3)$, $l \in \mathbb{R}$, we have the
Hamiltonian \\ $H(A,c,\Pi,\Gamma,\alpha,l): T^* Q \cong
\textmd{SE}(3)\times \mathfrak{se}^\ast
(3)\times\mathbb{R}\times\mathbb{R}\to \mathbb{R}$ given by
\begin{align}
 H(A,c,\Pi,\Gamma,\alpha,l) &=\Omega\cdot \Pi+\dot{\alpha}\cdot
l-L(A,c,\Omega,\Gamma,\alpha,\dot{\alpha}) \nonumber \\
&=\frac{1}{2}[\frac{\Pi_1^2}{\bar{I}_1}+\frac{\Pi_2^2}{\bar{I}_2}
+\frac{(\Pi_3-l)^2}{\bar{I}_3}+\frac{l^2}{J_3}]+ mgh\Gamma \cdot
\chi \; . \end{align}

In order to give the motion equation of spacecraft-rotor system, in
the following we need to consider the symmetry and reduced
symplectic structure of the configuration space $Q=
\textmd{SE}(3)\times S^1$.

\section{Symmetric Reduction of the Rigid Spacecraft with a Rotor}

In the following we regard the rigid spacecraft with an internal
rotor as a regular point reducible RCH system on the generalization
of rotation group $\textmd{SO}(3)\times S^1$ and on the
generalization of Euclidean group $\textmd{SE}(3)\times S^1$,
respectively, and give the motion equations of their reduced RCH
systems on the symplectic leaves by calculation in detail. It is
worthy note that they are different from the symmetric reductions of
Hamiltonian systems in Bloch and Leonard \cite{blle02}, Bloch,
Leonard and Marsden \cite{bllema01}, Marsden \cite{ma92}, the
reductions in this paper are all the controlled Hamiltonian
reductions, that is, the symmetric reductions of (regular)
controlled Hamiltonian systems, see Marsden et al \cite{mawazh10}.
We also follow the notations and conventions introduced in Marsden
et al \cite{mamiorpera07, mamora90}, Marsden and Ratiu
\cite{mara99}, Libermann and Marle \cite{lima87}, Ortega and Ratiu
\cite{orra04}.

\subsection{Symmetric Reduction of Spacecraft-Rotor System with Coincident Centers}

We first give the regular point reduction of spacecraft-rotor system
with coincident centers of buoyancy and gravity. Assume that Lie
group $G=\textmd{SO}(3)$ acts freely and properly on
$Q=\textmd{SO}(3)\times S^1$ by the left translation on the first
factor $\textmd{SO}(3)$, and the trivial action on the second factor
$S^1$. By using the left trivialization of $T^\ast \textmd{SO}(3)=
\textmd{SO}(3)\times \mathfrak{so}^\ast(3) $, then the action of
$\textmd{SO}(3)$ on phase space $T^\ast Q= T^\ast
\textmd{SO}(3)\times T^\ast S^1$ is by cotangent lift of left
translation on $\textmd{SO}(3)$ at the identity, that is, $\Phi:
\textmd{SO}(3)\times T^\ast \textmd{SO}(3)\times T^\ast S^1 \cong
\textmd{SO}(3)\times \textmd{SO}(3)\times \mathfrak{so}^\ast(3)
\times \mathbb{R} \times \mathbb{R}\to \textmd{SO}(3)\times
\mathfrak{so}^\ast(3)\times \mathbb{R} \times \mathbb{R},$ given by
$\Phi(B,(A,\Pi,\alpha,l))=(BA,\Pi,\alpha,l)$, for any $A,B\in
\textmd{SO}(3), \; \Pi \in \mathfrak{so}^\ast(3), \; \alpha,l \in
\mathbb{R}$, which is also free and proper, and the orbit space
$(T^* Q)/ \textmd{SO}(3)$ is a smooth manifold and $\pi: T^*Q
\rightarrow (T^*Q )/ \textmd{SO}(3) $ is a smooth submersion. Since
$\textmd{SO}(3)$ acts trivially on $\mathfrak{so}^\ast(3)$ and
$\mathbb{R}$, it follows that $(T^\ast Q)/ \textmd{SO}(3)$ is
diffeomorphic to $\mathfrak{so}^\ast(3) \times \mathbb{R}
\times \mathbb{R}$.\\

We know that $\mathfrak{so}^\ast(3)$ is a Poisson manifold with
respect to its rigid body Lie-Poisson bracket defined by
\begin{equation}
\{F,K\}_{\mathfrak{so}^\ast(3)}(\Pi)= -\Pi\cdot (\nabla_\Pi
F\times \nabla_\Pi K), \;\; \forall F,K\in
C^\infty(\mathfrak{so}^\ast(3)),\;\; \Pi \in
\mathfrak{so}^\ast(3).\label{3.1}
\end{equation}
For $\mu \in \mathfrak{so}^\ast(3)$, the coadjoint orbit
$\mathcal{O}_\mu \subset \mathfrak{so}^\ast(3)$ has the induced
orbit symplectic form $\omega^{-}_{\mathcal{O}_\mu}$, which is
coincide with the restriction of the Lie-Poisson bracket on
$\mathfrak{so}^\ast(3)$ to the coadjoint orbit $\mathcal{O}_\mu$.
From the Symplectic Stratification theorem we know that the
coadjoint orbits $(\mathcal{O}_\mu, \omega_{\mathcal{O}_\mu}^{-}),
\; \mu\in \mathfrak{so}^\ast(3),$ form the symplectic leaves of the
Poisson manifold
$(\mathfrak{so}^\ast(3),\{\cdot,\cdot\}_{\mathfrak{so}^\ast(3)}). $
Let $\omega_{\mathbb{R}}$ be the canonical symplectic form on
$T^\ast \mathbb{R} \cong \mathbb{R} \times \mathbb{R}$, which is
given by
\begin{equation}
 \omega_{\mathbb{R}}((\theta_1, \lambda_1),(\theta_2,
\lambda_2))=<\lambda_2,\theta_1> -<\lambda_1,\theta_2>,
\end{equation} where $(\theta_i, \lambda_i)\in \mathbb{R}\times
\mathbb{R}, \; i=1,2$, $<\cdot,\cdot>$ is the standard inner product
on $\mathbb{R}\times \mathbb{R}$. It induces a canonical Poisson
bracket $\{\cdot,\cdot\}_{\mathbb{R}}$ on $T^\ast \mathbb{R}$, which
is given by
\begin{equation}
\{F,K\}_{\mathbb{R}}(\theta,\lambda)= \frac{\partial
F}{\partial \theta} \frac{\partial K}{\partial \lambda}-
\frac{\partial K}{\partial \theta}\frac{\partial F}{\partial
\lambda} .\label{3.3}
\end{equation}
Thus, we can induce a symplectic form
$\tilde{\omega}^{-}_{\mathcal{O}_\mu \times \mathbb{R} \times
\mathbb{R}}= \pi_{\mathcal{O}_\mu}^\ast
\omega^{-}_{\mathcal{O}_\mu}+ \pi_{\mathbb{R}}^\ast
\omega_{\mathbb{R}}$ on the smooth manifold $\mathcal{O}_\mu \times
\mathbb{R} \times \mathbb{R}$, where the maps
$\pi_{\mathcal{O}_\mu}: \mathcal{O}_\mu \times \mathbb{R} \times
\mathbb{R} \to \mathcal{O}_\mu$ and $\pi_{\mathbb{R}}:
\mathcal{O}_\mu \times \mathbb{R} \times \mathbb{R} \to \mathbb{R}
\times \mathbb{R}$ are canonical projections, and induce a Poisson
bracket $\{\cdot,\cdot\}_{-}= \pi_{\mathfrak{so}^\ast(3)}^\ast
\{\cdot,\cdot\}_{\mathfrak{so}^\ast(3)}+ \pi_{\mathbb{R}}^\ast
\{\cdot,\cdot\}_{\mathbb{R}}$ on the smooth manifold
$\mathfrak{so}^\ast(3)\times \mathbb{R} \times \mathbb{R}$, where
the maps $\pi_{\mathfrak{so}^\ast(3)}: \mathfrak{so}^\ast(3) \times
\mathbb{R} \times \mathbb{R} \to \mathfrak{so}^\ast(3)$ and
$\pi_{\mathbb{R}}: \mathfrak{so}^\ast(3) \times \mathbb{R} \times
\mathbb{R} \to \mathbb{R} \times \mathbb{R}$ are canonical
projections, and such that $(\mathcal{O}_\mu \times \mathbb{R}\times
\mathbb{R},\tilde{\omega}_{\mathcal{O}_\mu \times \mathbb{R} \times
\mathbb{R}}^{-})$ is a symplectic leaf of the Poisson manifold
$(\mathfrak{so}^\ast(3) \times \mathbb{R} \times \mathbb{R}, \{\cdot,\cdot\}_{-}). $\\

On the other hand, from $T^\ast Q = T^\ast \textmd{SO}(3) \times
T^\ast S^1$ we know that there is a canonical symplectic form
$\omega_Q= \pi^\ast_{\textmd{SO}(3)} \omega_0 +\pi^\ast_{S^1}
\omega_{S^1}$ on $T^\ast Q$, where $\omega_0$ is the canonical
symplectic form on $T^\ast \textmd{SO}(3)$ and the maps
$\pi_{\textmd{SO}(3)}: Q= \textmd{SO}(3)\times S^1 \to
\textmd{SO}(3)$ and $\pi_{S^1}: Q= \textmd{SO}(3)\times S^1 \to S^1$
are canonical projections. Then the cotangent lift of left
$\textmd{SO}(3)$-action $\Phi: \textmd{SO}(3) \times T^\ast Q \to
T^\ast Q$ is also symplectic, and admits an associated
$\operatorname{Ad}^\ast$-equivariant momentum map $\mathbf{J}_Q:
T^\ast Q \to \mathfrak{so}^\ast(3)$ such that $\mathbf{J}_Q\cdot
\pi^\ast_{\textmd{SO}(3)}=\mathbf{J}_{\textmd{SO}(3)}$, where
$\mathbf{J}_{\textmd{SO}(3)}: T^\ast \textmd{SO}(3) \rightarrow
\mathfrak{so}^\ast(3)$ is a momentum map of left
$\textmd{SO}(3)$-action on $T^\ast \textmd{SO}(3)$, and
$\pi^\ast_{\textmd{SO}(3)}: T^\ast \textmd{SO}(3) \to T^\ast Q$. If
$\mu\in\mathfrak{so}^\ast(3)$ is a regular value of $\mathbf{J}_Q$,
then $\mu\in\mathfrak{so}^\ast(3)$ is also a regular value of
$\mathbf{J}_{\textmd{SO}(3)}$ and $\mathbf{J}_Q^{-1}(\mu)\cong
\mathbf{J}_{\textmd{SO}(3)}^{-1}(\mu)\times \mathbb{R} \times
\mathbb{R}$. Denote by $\textmd{SO}(3)_\mu=\{g\in
\textmd{SO}(3)|\operatorname{Ad}_g^\ast \mu=\mu \}$ the isotropy
subgroup of coadjoint $\textmd{SO}(3)$-action at the point
$\mu\in\mathfrak{so}^\ast(3)$. It follows that $\textmd{SO}(3)_\mu$
acts also freely and properly on $\mathbf{J}_Q^{-1}(\mu)$, the
regular point reduced space $(T^\ast
Q)_\mu=\mathbf{J}_Q^{-1}(\mu)/\textmd{SO}(3)_\mu\cong (T^\ast
\textmd{SO}(3))_\mu \times \mathbb{R} \times \mathbb{R}$ of $(T^\ast
Q,\omega_Q)$ at $\mu$, is a symplectic manifold with symplectic form
$\omega_\mu$ uniquely characterized by the relation $\pi_\mu^\ast
\omega_\mu=i_\mu^\ast \omega_Q=i_\mu^\ast \pi^\ast_{\textmd{SO}(3)}
\omega_0 +i_\mu^\ast \pi^\ast_{S^1} \omega_{S^1}$, where the map
$i_\mu:\mathbf{J}_Q^{-1}(\mu)\rightarrow T^\ast Q$ is the inclusion
and $\pi_\mu:\mathbf{J}_Q^{-1}(\mu)\rightarrow (T^\ast Q)_\mu$ is
the projection. Because from Abraham and Marsden \cite{abma78}, we
know that $((T^\ast \textmd{SO}(3))_\mu,\omega_\mu)$ is
symplectically diffeomorphic to
$(\mathcal{O}_\mu,\omega_{\mathcal{O}_\mu}^{-})$, and hence we have
that $((T^\ast Q)_\mu,\omega_\mu)$ is symplectically diffeomorphic
to $(\mathcal{O}_\mu \times \mathbb{R}\times
\mathbb{R},\tilde{\omega}_{\mathcal{O}_\mu \times \mathbb{R} \times
\mathbb{R}}^{-})$, which is a symplectic leaf of the Poisson
manifold
$(\mathfrak{so}^\ast(3) \times \mathbb{R} \times \mathbb{R}, \{\cdot,\cdot\}_{-}). $\\

From the expression $(2.1)$ of the Hamiltonian, we know that
$H(A,\Pi,\alpha,l)$ is invariant under the left
$\textmd{SO}(3)$-action $\Phi: \textmd{SO}(3)\times T^\ast Q \to
T^\ast Q$. For the case $\mu \in \mathfrak{so}^\ast(3)$ is the
regular value of $\mathbf{J}_Q$, we have the reduced Hamiltonian
$h_\mu(\Pi, \alpha,l):\mathcal{O}_\mu
\times\mathbb{R}\times\mathbb{R}(\subset \mathfrak{so}^\ast
(3)\times\mathbb{R}\times\mathbb{R})\to \mathbb{R}$ given by
$h_\mu(\Pi,\alpha,l)= \pi_\mu(H(A,\Pi,\alpha,l))
=H(A,\Pi,\alpha,l)|_{\mathcal{O}_\mu
\times\mathbb{R}\times\mathbb{R}}.$ From the rigid body Poisson
bracket on $\mathfrak{so}^\ast(3)$ and the Poisson bracket on
$T^\ast \mathbb{R}$, we can get the Poisson bracket on
$\mathfrak{so}^\ast(3)\times\mathbb{R}\times\mathbb{R}$, that is,
for $F,K: \mathfrak{so}^\ast(3)\times\mathbb{R}\times\mathbb{R} \to
\mathbb{R}, $ we have that
\begin{align}
\{F,K\}_{-}(\Pi,\alpha,l)=-\Pi\cdot(\nabla_\Pi F\times \nabla_\Pi
K)+ \{F,K\}_{\mathbb{R}}(\alpha,l).
\end{align}
See Krishnaprasad and Marsden \cite{krma87}. In particular, for
$F_\mu,K_\mu: \mathcal{O}_\mu \times\mathbb{R}\times\mathbb{R} \to
\mathbb{R}$, we have that $\tilde{\omega}_{\mathcal{O}_\mu
\times\mathbb{R}\times\mathbb{R}}^{-}(X_{F_\mu}, X_{K_\mu})=
\{F_\mu,K_\mu\}_{-}|_{\mathcal{O}_\mu
\times\mathbb{R}\times\mathbb{R} }$. Moreover, for reduced
Hamiltonian $h_\mu(\Pi,\alpha,l): \mathcal{O}_\mu
\times\mathbb{R}\times\mathbb{R} \to \mathbb{R}$, we have the
Hamiltonian vector field
$X_{h_\mu}(K_\mu)=\{K_\mu,h_\mu\}_{-}|_{\mathcal{O}_\mu
\times\mathbb{R}\times\mathbb{R} },$ and hence we have that
\begin{align*}
& X_{h_\mu}(\Pi)(\Pi,\alpha,l)=\{\Pi,h_\mu\}_{-}(\Pi,\alpha,l)\\
& =-\Pi\cdot(\nabla_\Pi\Pi\times\nabla_\Pi h_\mu)+ (\frac{\partial
\Pi}{\partial \alpha} \frac{\partial h_\mu}{\partial l}-
\frac{\partial h_\mu}{\partial \alpha}\frac{\partial \Pi}{\partial
l}) \\
& =-\nabla_\Pi\Pi\cdot(\nabla_\Pi h_\mu\times
\Pi)=\Pi\times\Omega,
\end{align*}
since $\nabla_\Pi\Pi=1, \; \frac{\partial \Pi}{\partial \alpha}=
\frac{\partial \Pi}{\partial l}=0, $ and $\nabla_\Pi h_{\mu}=\Omega.
$

\begin{align*}
& X_{h_\mu}(\alpha)(\Pi,\alpha,l)=\{\alpha,h_\mu\}_{-}(\Pi,\alpha,l)\\
& =-\Pi\cdot(\nabla_\Pi\alpha \times \nabla_\Pi h_\mu)+
(\frac{\partial \alpha}{\partial \alpha} \frac{\partial
h_\mu}{\partial l}- \frac{\partial h_\mu}{\partial
\alpha}\frac{\partial \alpha}{\partial l}) \\
& = -\frac{(\Pi_3-l)}{\bar{I}_3}+\frac{l}{J_3},
\end{align*}
since $\nabla_\Pi\alpha=0, \; \frac{\partial \alpha}{\partial l}=0,
$ and $\frac{\partial h_{\mu}}{\partial l}= -\frac{(\Pi_3-
l)}{\bar{I}_3} +\frac{l}{J_3}$.

\begin{align*}
& X_{h_\mu}(l)(\Pi,\alpha,l)=\{l,h_\mu\}_{-}(\Pi,\alpha,l)\\
& =-\Pi\cdot(\nabla_\Pi l \times \nabla_\Pi h_\mu)+ (\frac{\partial
l}{\partial \alpha} \frac{\partial h_\mu}{\partial l}-
\frac{\partial h_\mu}{\partial \alpha}\frac{\partial l}{\partial l})
= 0,
\end{align*}
since $\nabla_\Pi l=0, $ and $\frac{\partial l}{\partial \alpha}=
\frac{\partial h_\mu}{\partial \alpha}=0. $ If we consider the rigid
spacecraft-rotor system with a control torque $u: T^\ast Q \to W $
acting on the rotor, and $u\in W \subset \mathbf{J}^{-1}_Q(\mu)$ is
invariant under the left $\textmd{SO}(3)$-action, and its reduced
control torque $u_\mu: \mathcal{O}_\mu
\times\mathbb{R}\times\mathbb{R} \to W_\mu (\subset \mathcal{O}_\mu
\times\mathbb{R}\times\mathbb{R}) $ is given by
$u_\mu(\Pi,\alpha,l)= \pi_\mu(u(A,\Pi,\alpha,l))=
u(A,\Pi,\alpha,l)|_{\mathcal{O}_\mu \times\mathbb{R}\times\mathbb{R}
}, $ where $\pi_\mu: \mathbf{J}_Q^{-1}(\mu) \rightarrow
\mathcal{O}_\mu \times\mathbb{R}\times\mathbb{R},\; W_\mu
=\pi_\mu(W). $ Moreover, the dynamical vector field of regular point
reduced spacecraft-rotor system $(\mathcal{O}_\mu \times \mathbb{R}
\times \mathbb{R},\tilde{\omega}_{\mathcal{O}_\mu \times \mathbb{R}
\times \mathbb{R}}^{-},h_\mu,u_\mu)$ is given by
$$ X_{(\mathcal{O}_\mu \times \mathbb{R} \times
\mathbb{R},\tilde{\omega}_{\mathcal{O}_\mu \times \mathbb{R} \times
\mathbb{R}}^{-},h_\mu,u_\mu)} = X_{h_\mu} + \mbox{vlift}(u_\mu) ,
$$
where $\mbox{vlift}(u_\mu)= \mbox{vlift}(u_\mu)X_{h_\mu} \in
T(\mathcal{O}_\mu \times\mathbb{R}\times\mathbb{R}). $ Note that
$\mbox{vlift}(u_\mu)X_{h_\mu}$ is the vertical lift of vector field
$X_{h_\mu}$ under the action of $u_\mu$ along fibers, that is,
$$\mbox{vlift}(u_\mu)X_{h_\mu}(\Pi,\alpha,l)
= \mbox{vlift}((Tu_\mu X_{h_\mu})(u_\mu(\Pi,\alpha,l)),
(\Pi,\alpha,l)) = (Tu_\mu X_{h_\mu})^v_\sigma(\Pi,\alpha,l),$$ where
$\sigma$ is a geodesic in $\mathcal{O}_\mu
\times\mathbb{R}\times\mathbb{R} $ connecting $u_\mu(\Pi,\alpha,l)$
and $(\Pi,\alpha,l)$, and $(Tu_\mu
X_{h_\mu})^v_\sigma(\Pi,\alpha,l)$ is the parallel displacement of
the vertical vector $(Tu_\mu X_{h_\mu})^v(\Pi,\alpha,l)$ along the
geodesic $\sigma$ from $u_\mu(\Pi,\alpha,l)$ to $(\Pi,\alpha,l)$,
see Marsden et al \cite{mawazh10} and Wang \cite{wa13d}. If we
assume that
$$\mbox{vlift}(u_{\mu})X_{h_{\mu}}(\Pi, \alpha, l)
 =(\mathcal{U}_\Pi, \mathcal{U}_\alpha, \mathcal{U}_l) \in T_x(\mathcal{O}_{\mu} \times\mathbb{R}\times\mathbb{R}),
$$ where $x= (\Pi, \alpha, l) \in \mathcal{O}_{\mu} \times\mathbb{R}\times\mathbb{R}, $
and $\mathcal{U}_\Pi \in \mathbb{R}^3, $ and $\mathcal{U}_\alpha, \;
\mathcal{U}_l \in \mathbb{R}. $ Thus, in the case of coincident
centers of buoyancy and gravity, the equations of motion for reduced
spacecraft-rotor system with the control torque $u$ acting on the
rotor are given by
\begin{equation}
\left\{\begin{aligned} & \frac{\mathrm{d}\Pi}{\mathrm{d}t}=\Pi\times\Omega + \mathcal{U}_\Pi, \\
& \frac{\mathrm{d}\alpha}{\mathrm{d}t} =
-\frac{(\Pi_3-l)}{\bar{I}_3}+\frac{l}{J_3} + \mathcal{U}_\alpha, \\
& \frac{\mathrm{d}l}{\mathrm{d}t}= \mathcal{U}_l,
\end{aligned}\right.\label{3.5}\end{equation}

To sum up the above discussion, we have the following theorem.
\begin{theo}
In the case of coincident centers of buoyancy and gravity, the
spacecraft-rotor system with the control torque $u$ acting on the
rotor, that is, the 5-tuple $(T^\ast Q, \textmd{SO}(3), \omega_Q, H,
u ), $ where $Q= \textmd{SO}(3)\times S^1, $ is a regular point
reducible RCH system. For a point $\mu \in \mathfrak{so}^\ast(3)$,
the regular value of the momentum map $\mathbf{J}_Q:
\textmd{SO}(3)\times \mathfrak{so}^\ast(3) \times \mathbb{R} \times
\mathbb{R} \to \mathfrak{so}^\ast(3)$, the regular point reduced
system is the 4-tuple $(\mathcal{O}_\mu \times \mathbb{R} \times
\mathbb{R},\tilde{\omega}_{\mathcal{O}_\mu \times \mathbb{R} \times
\mathbb{R}}^{-},h_\mu,u_\mu), $ where $\mathcal{O}_\mu \subset
\mathfrak{so}^\ast(3)$ is the coadjoint orbit,
$\tilde{\omega}_{\mathcal{O}_\mu \times \mathbb{R} \times
\mathbb{R}}^{-}$ is orbit symplectic form on $\mathcal{O}_\mu \times
\mathbb{R}\times \mathbb{R} $,
$h_\mu(\Pi,\alpha,l)=\pi_\mu(H(A,\Pi,\alpha,l))=H(A,\Pi,\alpha,l)|_{\mathcal{O}_\mu
\times\mathbb{R}\times\mathbb{R}}$, $u_\mu(\Pi,\alpha,l)=
\pi_\mu(u(A,\Pi,\alpha,l))= u(A,\Pi,\alpha,l)|_{\mathcal{O}_\mu
\times\mathbb{R}\times\mathbb{R}}$, and its equations of motion are
given by (\ref{3.5}).
\end{theo}

\begin{rema}
When the rigid spacecraft does not carry any internal rotor, in this
case the configuration space is $Q=G=\textmd{SO}(3), $ the motion of
rigid spacecraft is just the rotation motion of a rigid body, the
above symmetric reduction of spacecraft-rotor system is just the
Marsden-Weinstein reduction of a rigid body at a regular value of
momentum map, and the motion equation $(3.5)$ of reduced
spacecraft-rotor system becomes the motion equation of a reduced
rigid body on a coadjoint orbit of Lie group $\textmd{SO}(3)$. See
Marsden and Ratiu \cite{mara99}.
\end{rema}

\subsection{Symmetric Reduction of Spacecraft-Rotor System with Non-coincident Centers}

In the following we shall give the regular point reduction of
spacecraft-rotor system with non-coincident centers of buoyancy and
gravity. Because the drift in the direction of gravity breaks the
symmetry and the spacecraft-rotor system is no longer
$\textmd{SO}(3)$ invariant. In this case, its physical phase space
is $T^\ast \textmd{SO}(3)\times T^* S^1$ and the symmetry group is
$S^1$, regarded as rotations about the third principal axis, that
is, the axis of gravity. By the semidirect product reduction
theorem, see Marsden et al \cite{mamiorpera07}, we know that the
reduction of $T^\ast \textmd{SO}(3)$ by $S^1$ gives a space which is
symplectically diffeomorphic to the reduced space obtained by the
reduction of $T^\ast \textmd{SE}(3)$ by left action of
$\textmd{SE}(3)$, that is the coadjoint orbit $\mathcal{O}_{(\mu,a)}
\subset \mathfrak{se}^\ast(3)\cong T^\ast
\textmd{SE}(3)/\textmd{SE}(3)$. In fact, in this case, we can
identify the phase space $T^\ast \textmd{SO}(3)$ with the reduction
of the cotangent bundle of the special Euclidean group
$\textmd{SE}(3)=\textmd{SO}(3)\circledS \mathbb{R}^3$ by the
Euclidean translation subgroup $\mathbb{R}^3$ and identifies the
symmetry group $S^1$ with isotropy group $G_a=\{ A\in
\textmd{SO}(3)\mid Aa=a \}=S^1$, which is Abelian and
$(G_a)_{\mu_a}= G_a =S^1,\; \forall \mu_a \in \mathfrak{g}^\ast_a$,
where $a$ is a vector aligned with the direction of gravity and
where $\textmd{SO}(3)$ acts on $\mathbb{R}^3$ in the standard way.\\

Assume that Lie group $G=\textmd{SE}(3)$ acts freely and properly on
$Q=\textmd{SE}(3)\times S^1$ by the left translation on the first
factor $\textmd{SE}(3)$, and the trivial action on the second factor
$S^1$. By using the left trivialization of $T^\ast \textmd{SE}(3)=
\textmd{SE}(3)\times \mathfrak{se}^\ast(3) $, then the action of
$\textmd{SE}(3)$ on phase space $T^\ast Q= T^\ast
\textmd{SE}(3)\times T^\ast S^1$ is by cotangent lift of left
translation on $\textmd{SE}(3)$ at the identity, that is, $\Phi:
\textmd{SE}(3)\times T^\ast \textmd{SE}(3)\times T^\ast S^1 \cong
\textmd{SE}(3)\times \textmd{SE}(3)\times \mathfrak{se}^\ast(3)
\times \mathbb{R} \times \mathbb{R}\to \textmd{SE}(3)\times
\mathfrak{se}^\ast(3)\times \mathbb{R} \times \mathbb{R},$ given by
$\Phi((B,b)(A,c,\Pi,\Gamma,\alpha,l))=(BA,c,\Pi,\Gamma,\alpha,l)$,
for any $A,B\in \textmd{SO}(3), \; \Pi \in \mathfrak{so}^\ast(3), \;
b, c, \Gamma \in \mathbb{R}^3, \; \alpha,l \in \mathbb{R}$, which is
also free and proper, and the orbit space $(T^* Q)/ \textmd{SE}(3)$
is a smooth manifold and $\pi: T^*Q \rightarrow (T^*Q )/
\textmd{SE}(3) $ is a smooth submersion. Since $\textmd{SE}(3)$ acts
trivially on $\mathfrak{se}^\ast(3)$ and $\mathbb{R}$, it follows
that $(T^\ast Q)/ \textmd{SE}(3)$ is diffeomorphic to
$\mathfrak{se}^\ast(3) \times \mathbb{R}
\times \mathbb{R}$.\\

We know that $\mathfrak{se}^\ast(3)$ is a Poisson manifold with
respect to its heavy top Lie-Poisson bracket defined by
\begin{equation}   \{F,K\}_{\mathfrak{se}^\ast(3)}(\Pi,\Gamma)= -\Pi\cdot (\nabla_\Pi
F\times \nabla_\Pi K)-\Gamma\cdot(\nabla_\Pi F\times \nabla_\Gamma
K-\nabla_\Pi K\times \nabla_\Gamma F), \;\; \label{3.6}
\end{equation}
where $F, K \in C^\infty(\mathfrak{se}^\ast(3)), \; (\Pi,\Gamma) \in
\mathfrak{se}^\ast(3)$. For $(\mu,a) \in \mathfrak{se}^\ast(3)$, the
coadjoint orbit $\mathcal{O}_{(\mu,a)} \subset
\mathfrak{se}^\ast(3)$ has the induced orbit symplectic form
$\omega^{-}_{\mathcal{O}_{(\mu,a)}}$, which is coincide with the
restriction of the Lie-Poisson bracket on $\mathfrak{se}^\ast(3)$ to
the coadjoint orbit $\mathcal{O}_{(\mu,a)}$, and the coadjoint
orbits $(\mathcal{O}_{(\mu,a)}, \omega_{\mathcal{O}_{(\mu,a)}}^{-}),
\; (\mu,a)\in \mathfrak{se}^\ast(3),$ form the symplectic leaves of
the Poisson manifold\\
$(\mathfrak{se}^\ast(3),\{\cdot,\cdot\}_{\mathfrak{se}^\ast(3)}). $
Let $\omega_{\mathbb{R}}$ be the canonical symplectic form on
$T^\ast \mathbb{R} \cong \mathbb{R} \times \mathbb{R}$ given by
$(3.2)$, and it induces a canonical Poisson bracket
$\{\cdot,\cdot\}_{\mathbb{R}}$ on $T^\ast \mathbb{R}$ given by
$(3.3)$. Thus, we can induce a symplectic form
$\tilde{\omega}^{-}_{\mathcal{O}_{(\mu,a)} \times \mathbb{R} \times
\mathbb{R}}= \pi_{\mathcal{O}_{(\mu,a)}}^\ast
\omega^{-}_{\mathcal{O}_{(\mu,a)}}+ \pi_{\mathbb{R}}^\ast
\omega_{\mathbb{R}}$ on the smooth manifold $\mathcal{O}_{(\mu,a)}
\times \mathbb{R} \times \mathbb{R}$, where the maps
$\pi_{\mathcal{O}_{(\mu,a)}}: \mathcal{O}_{(\mu,a)} \times
\mathbb{R} \times \mathbb{R} \to \mathcal{O}_{(\mu,a)}$ and
$\pi_{\mathbb{R}}: \mathcal{O}_{(\mu,a)} \times \mathbb{R} \times
\mathbb{R} \to \mathbb{R} \times \mathbb{R}$ are canonical
projections, and induce a Poisson bracket $\{\cdot,\cdot\}_{-}=
\pi_{\mathfrak{se}^\ast(3)}^\ast
\{\cdot,\cdot\}_{\mathfrak{se}^\ast(3)}+ \pi_{\mathbb{R}}^\ast
\{\cdot,\cdot\}_{\mathbb{R}}$ on the smooth manifold
$\mathfrak{se}^\ast(3)\times \mathbb{R} \times \mathbb{R}$, where
the maps $\pi_{\mathfrak{se}^\ast(3)}: \mathfrak{se}^\ast(3) \times
\mathbb{R} \times \mathbb{R} \to \mathfrak{se}^\ast(3)$ and
$\pi_{\mathbb{R}}: \mathfrak{se}^\ast(3) \times \mathbb{R} \times
\mathbb{R} \to \mathbb{R} \times \mathbb{R}$ are canonical
projections, and such that $(\mathcal{O}_{(\mu,a)} \times
\mathbb{R}\times \mathbb{R},\tilde{\omega}_{\mathcal{O}_{(\mu,a)}
\times \mathbb{R} \times \mathbb{R}}^{-})$ is a symplectic leaf of
the Poisson manifold
$(\mathfrak{se}^\ast(3) \times \mathbb{R} \times \mathbb{R}, \{\cdot,\cdot\}_{-}). $\\

On the other hand, from $T^\ast Q = T^\ast \textmd{SE}(3) \times
T^\ast S^1$ we know that there is a canonical symplectic form
$\omega_Q= \pi^\ast_{\textmd{SE}(3)} \omega_1 +\pi^\ast_{S^1}
\omega_{S^1}$ on $T^\ast Q$, where $\omega_1$ is the canonical
symplectic form on $T^\ast \textmd{SE}(3)$ and the maps
$\pi_{\textmd{SE}(3)}: Q= \textmd{SE}(3)\times S^1 \to
\textmd{SE}(3)$ and $\pi_{S^1}: Q= \textmd{SE}(3)\times S^1 \to S^1$
are canonical projections. Then the cotangent lift of left
$\textmd{SE}(3)$-action $\Phi: \textmd{SE}(3) \times T^\ast Q \to
T^\ast Q$ is also symplectic, and admits an associated
$\operatorname{Ad}^\ast$-equivariant momentum map $\mathbf{J}_Q:
T^\ast Q \to \mathfrak{se}^\ast(3)$ such that $\mathbf{J}_Q\cdot
\pi^\ast_{\textmd{SE}(3)}=\mathbf{J}_{\textmd{SE}(3)}$, where
$\mathbf{J}_{\textmd{SE}(3)}: T^\ast \textmd{SE}(3) \rightarrow
\mathfrak{se}^\ast(3)$ is a momentum map of left
$\textmd{SE}(3)$-action on $T^\ast \textmd{SE}(3)$, and
$\pi^\ast_{\textmd{SE}(3)}: T^\ast \textmd{SE}(3) \to T^\ast Q$. If
$(\mu,a)\in\mathfrak{se}^\ast(3)$ is a regular value of
$\mathbf{J}_Q$, then $(\mu,a)\in\mathfrak{se}^\ast(3)$ is also a
regular value of $\mathbf{J}_{\textmd{SE}(3)}$ and
$\mathbf{J}_Q^{-1}(\mu,a)\cong
\mathbf{J}_{\textmd{SE}(3)}^{-1}(\mu,a)\times \mathbb{R} \times
\mathbb{R}$. Denote by $\textmd{SE}(3)_{(\mu,a)}=\{g\in
\textmd{SE}(3)|\operatorname{Ad}_g^\ast (\mu,a)=(\mu,a) \}$ the
isotropy subgroup of coadjoint $\textmd{SE}(3)$-action at the point
$(\mu,a)\in\mathfrak{se}^\ast(3)$. It follows that
$\textmd{SE}(3)_{(\mu,a)}$ acts also freely and properly on
$\mathbf{J}_Q^{-1}(\mu,a)$, the regular point reduced space $(T^\ast
Q)_{(\mu,a)}=\mathbf{J}_Q^{-1}(\mu,a)/\textmd{SE}(3)_{(\mu,a)}\cong
(T^\ast \textmd{SE}(3))_{(\mu,a)} \times \mathbb{R} \times
\mathbb{R}$ of $(T^\ast Q,\omega_Q)$ at $(\mu,a)$, is a symplectic
manifold with symplectic form $\omega_{(\mu,a)}$ uniquely
characterized by the relation $\pi_{(\mu,a)}^\ast
\omega_{(\mu,a)}=i_{(\mu,a)}^\ast \omega_Q=i_{(\mu,a)}^\ast
\pi^\ast_{\textmd{SE}(3)} \omega_1 +i_{(\mu,a)}^\ast \pi^\ast_{S^1}
\omega_{S^1}$, where the map
$i_{(\mu,a)}:\mathbf{J}_Q^{-1}(\mu,a)\rightarrow T^\ast Q$ is the
inclusion and $\pi_{(\mu,a)}:\mathbf{J}_Q^{-1}(\mu,a)\rightarrow
(T^\ast Q)_{(\mu,a)}$ is the projection. Because from Abraham and
Marsden \cite{abma78}, we know that $((T^\ast
\textmd{SE}(3))_{(\mu,a)},\omega_{(\mu,a)})$ is symplectically
diffeomorphic to
$(\mathcal{O}_{(\mu,a)},\omega_{\mathcal{O}_{(\mu,a)}}^{-})$, and
hence we have that $((T^\ast Q)_{(\mu,a)},\omega_{(\mu,a)})$ is
symplectically diffeomorphic to $(\mathcal{O}_{(\mu,a)} \times
\mathbb{R}\times \mathbb{R},\tilde{\omega}_{\mathcal{O}_{(\mu,a)}
\times \mathbb{R} \times \mathbb{R}}^{-})$, which is a symplectic
leaf of the Poisson manifold
$(\mathfrak{se}^\ast(3) \times \mathbb{R} \times \mathbb{R}, \{\cdot,\cdot\}_{-}). $\\

From the expression $(2.2)$ of the Hamiltonian, we know that
$H(A,c,\Pi,\Gamma,\alpha,l)$ is invariant under the left
$\textmd{SE}(3)$-action $\Phi: \textmd{SE}(3)\times T^\ast Q \to
T^\ast Q$. For the case $(\mu,a) \in \mathfrak{se}^\ast(3)$ is the
regular value of $\mathbf{J}_Q$, we have the reduced Hamiltonian
$h_{(\mu,a)}(\Pi,\Gamma, \alpha,l): \mathcal{O}_{(\mu,a)}
\times\mathbb{R}\times\mathbb{R}(\subset \mathfrak{se}^\ast
(3)\times\mathbb{R}\times\mathbb{R})\to \mathbb{R}$ given by
$h_{(\mu,a)}(\Pi,\Gamma,\alpha,l)=
\pi_{(\mu,a)}(H(A,c,\Pi,\Gamma,\alpha,l))=
H(A,c,\Pi,\Gamma,\alpha,l)|_{\mathcal{O}_{(\mu,a)}
\times\mathbb{R}\times\mathbb{R}}.$ From the heavy top Poisson
bracket on $\mathfrak{se}^\ast(3)$ and the Poisson bracket on
$T^\ast \mathbb{R}$, we can get the Poisson bracket on
$\mathfrak{se}^\ast(3)\times\mathbb{R}\times\mathbb{R}$, that is,
for $F,K: \mathfrak{se}^\ast(3)\times\mathbb{R}\times\mathbb{R} \to
\mathbb{R}, $ we have that
\begin{align}
\{F,K\}_{-}(\Pi,\Gamma,\alpha,l)= &-\Pi\cdot(\nabla_\Pi F\times
\nabla_\Pi K) \nonumber \\
&-\Gamma\cdot(\nabla_\Pi F\times
\nabla_\Gamma K-\nabla_\Pi K\times \nabla_\Gamma F)+
\{F,K\}_{\mathbb{R}}(\alpha,l).
\end{align}
See Krishnaprasad and Marsden \cite{krma87}. In particular, for
$F_{(\mu,a)},K_{(\mu,a)}: \mathcal{O}_{(\mu,a)}
\times\mathbb{R}\times\mathbb{R} \to \mathbb{R}$, we have that
$\tilde{\omega}_{\mathcal{O}_{(\mu,a)}
\times\mathbb{R}\times\mathbb{R}}^{-}(X_{F_{(\mu,a)}},
X_{K_{(\mu,a)}})=
\{F_{(\mu,a)},K_{(\mu,a)}\}_{-}|_{\mathcal{O}_{(\mu,a)}
\times\mathbb{R}\times\mathbb{R} }$. Moreover, for reduced
Hamiltonian $h_{(\mu,a)}(\Pi,\Gamma,\alpha,l): \mathcal{O}_{(\mu,a)}
\times\mathbb{R}\times\mathbb{R} \to \mathbb{R}$, we have the
Hamiltonian vector field
$X_{h_{(\mu,a)}}(K_{(\mu,a)})=\{K_{(\mu,a)},h_{(\mu,a)}\}_{-}|_{\mathcal{O}_{(\mu,a)}
\times\mathbb{R}\times\mathbb{R} },$ and hence we have that
\begin{align*}
& X_{h_{(\mu,a)}}(\Pi)(\Pi,\Gamma,\alpha,l)
=\{\Pi,h_{(\mu,a)}\}_{-}(\Pi,\Gamma,\alpha,l)\\
& =-\Pi\cdot(\nabla_\Pi\Pi \times \nabla_\Pi h_{(\mu,a)})
-\Gamma\cdot(\nabla_\Pi\Pi\times\nabla_\Gamma h_{(\mu,a)}-\nabla_\Pi
h_{(\mu,a)}\times\nabla_\Gamma\Pi)\\
& \;\;\;\; + (\frac{\partial \Pi}{\partial \alpha} \frac{\partial
h_{(\mu,a)}}{\partial l}- \frac{\partial h_{(\mu,a)}}{\partial
\alpha}\frac{\partial \Pi}{\partial l})\\ &
=\Pi\times\Omega-mgh\chi\times\Gamma=\Pi\times\Omega+mgh\Gamma\times\chi,
\end{align*}
since $\nabla_\Pi\Pi=1, \; \nabla_\Gamma\Pi =0, \; \frac{\partial
\Pi}{\partial \alpha}= \frac{\partial \Pi}{\partial l}=0, $ and
$\nabla_\Pi h_{(\mu,a)}=\Omega, \; \nabla_\Gamma h_{(\mu,a)}=
mgh\chi. $

\begin{align*}
& X_{h_{(\mu,a)}}(\Gamma)(\Pi,\Gamma,\alpha,l)
=\{\Gamma,h_{(\mu,a)}\}_{-}(\Pi,\Gamma,\alpha,l)\\
& =-\Pi\cdot(\nabla_\Pi\Gamma\times\nabla_\Pi h_{(\mu,a)})
-\Gamma\cdot(\nabla_\Pi\Gamma\times\nabla_\Gamma
h_{(\mu,a)}-\nabla_\Pi h_{(\mu,a)}\times\nabla_\Gamma\Gamma)\\
& \;\;\;\; + (\frac{\partial \Gamma}{\partial \alpha} \frac{\partial
h_{(\mu,a)}}{\partial l}- \frac{\partial h_{(\mu,a)}}{\partial
\alpha}\frac{\partial \Gamma}{\partial l})\\
& =\nabla_\Gamma\Gamma\cdot(\Gamma\times\nabla_\Pi
h_{(\mu,a)})=\Gamma\times\Omega,
\end{align*}
since $\nabla_\Gamma \Gamma =1, \; \nabla_\Pi\Gamma =0, \;
\frac{\partial \Gamma}{\partial \alpha}= \frac{\partial
\Gamma}{\partial l}=0. $

\begin{align*}
& X_{h_{(\mu,a)}}(\alpha)(\Pi,\Gamma,\alpha,l)
=\{\alpha,h_{(\mu,a)}\}_{-}(\Pi,\Gamma,\alpha,l)\\
& =-\Pi\cdot(\nabla_\Pi\alpha \times \nabla_\Pi h_{(\mu,a)})
-\Gamma\cdot(\nabla_\Pi \alpha \times\nabla_\Gamma
h_{(\mu,a)}-\nabla_\Pi h_{(\mu,a)}\times\nabla_\Gamma \alpha)\\
& \;\;\;\; +(\frac{\partial \alpha}{\partial \alpha} \frac{\partial
h_{(\mu,a)}}{\partial l}- \frac{\partial h_{(\mu,a)}}{\partial
\alpha}\frac{\partial \alpha}{\partial l}) \\
& = -\frac{(\Pi_3-l)}{\bar{I}_3}+\frac{l}{J_3},
\end{align*}
since $\nabla_\Pi\alpha= \nabla_\Gamma \alpha=0, \; \frac{\partial
\alpha}{\partial l}=0, $ and $\frac{\partial h_{(\mu,a)}}{\partial
l}= -\frac{(\Pi_3- l)}{\bar{I}_3} +\frac{l}{J_3}$.

\begin{align*}
& X_{h_{(\mu,a)}}(l)(\Pi,\Gamma,\alpha,l)
=\{l,h_{(\mu,a)}\}_{-}(\Pi,\Gamma,\alpha,l)\\
& =-\Pi\cdot(\nabla_\Pi l \times \nabla_\Pi h_{(\mu,a)})
-\Gamma\cdot(\nabla_\Pi l \times\nabla_\Gamma h_{(\mu,a)}-\nabla_\Pi
h_{(\mu,a)}\times\nabla_\Gamma l)\\
& \;\;\;\; + (\frac{\partial l}{\partial \alpha} \frac{\partial
h_{(\mu,a)}}{\partial l}- \frac{\partial h_{(\mu,a)}}{\partial
\alpha}\frac{\partial l}{\partial l})\\
& = 0,
\end{align*}
since $\nabla_\Pi l= \nabla_\Gamma l=0, $ and $\frac{\partial
l}{\partial \alpha}=\frac{\partial h_{(\mu,a)}}{\partial \alpha}=0.
$ If we consider the rigid spacecraft-rotor system with a control
torque $u: T^\ast Q \to W $ acting on the rotor, and $u\in W \subset
\mathbf{J}^{-1}_Q(\mu,a)$ is invariant under the left
$\textmd{SE}(3)$-action, and its reduced control torque
$u_{(\mu,a)}: \mathcal{O}_{(\mu,a)} \times\mathbb{R}\times\mathbb{R}
\to W_{(\mu,a)} (\subset \mathcal{O}_{(\mu,a)}
\times\mathbb{R}\times\mathbb{R}) $ is given by
$u_{(\mu,a)}(\Pi,\Gamma,\alpha,l)=
\pi_{(\mu,a)}(u(A,c,\Pi,\Gamma,\alpha,l))=
u(A,c,\Pi,\Gamma,\alpha,l)|_{\mathcal{O}_{(\mu,a)}
\times\mathbb{R}\times\mathbb{R} }, $ where $\pi_{(\mu,a)}:
\mathbf{J}_Q^{-1}(\mu,a) \rightarrow \mathcal{O}_{(\mu,a)}
\times\mathbb{R}\times\mathbb{R}, \; W_{(\mu,a)} = \pi_{(\mu,a)}(W).
$ Moreover, the dynamical vector field of regular point reduced
spacecraft-rotor system $(\mathcal{O}_{(\mu,a)} \times \mathbb{R}
\times \mathbb{R},\tilde{\omega}_{\mathcal{O}_{(\mu,a)} \times
\mathbb{R} \times \mathbb{R}}^{-},h_{(\mu,a)},u_{(\mu,a)})$ is given
by
$$ X_{(\mathcal{O}_{(\mu,a)} \times \mathbb{R} \times
\mathbb{R},\tilde{\omega}_{\mathcal{O}_{(\mu,a)} \times \mathbb{R}
\times \mathbb{R}}^{-},h_{(\mu,a)},u_{(\mu,a)})} = X_{h_{(\mu,a)}} +
\mbox{vlift}(u_{(\mu,a)}), $$ where $\mbox{vlift}(u_{(\mu,a)})=
\mbox{vlift}(u_{(\mu,a)})X_{h_{(\mu,a)}} \in T(\mathcal{O}_{(\mu,a)}
\times\mathbb{R}\times\mathbb{R}). $ Note that
$\mbox{vlift}(u_{(\mu,a)})X_{h_{(\mu,a)}}$ is the vertical lift of
vector field $X_{h_{(\mu,a)}}$ under the action of $u_{(\mu,a)}$
along fibers, that is,
\begin{align*}
\mbox{vlift}(u_{(\mu,a)})X_{h_{(\mu,a)}}(\Pi,\Gamma,\alpha,l) & =
\mbox{vlift}((Tu_{(\mu,a)}
X_{h_{(\mu,a)}})(u_{(\mu,a)}(\Pi,\Gamma,\alpha,l)),
(\Pi,\Gamma,\alpha,l))\\ & = (Tu_{(\mu,a)}
X_{h_{(\mu,a)}})^v_\sigma(\Pi,\Gamma,\alpha,l),
\end{align*}
where $\sigma$ is a geodesic in $\mathcal{O}_{(\mu,a)}
\times\mathbb{R}\times\mathbb{R} $ connecting
$u_{(\mu,a)}(\Pi,\Gamma,\alpha,l)$ and $(\Pi,\Gamma,\alpha,l)$, and
\\ $(Tu_{(\mu,a)} X_{h_{(\mu,a)}})^v_\sigma(\Pi,\Gamma,\alpha,l)$ is
the parallel displacement of the vertical vector \\ $(Tu_{(\mu,A)}
X_{h_{(\mu,a)}})^v(\Pi,\Gamma,\alpha,l)$ along the geodesic $\sigma$
from $u_{(\mu,a)}(\Pi,\Gamma,\alpha,l)$ to $(\Pi,\Gamma,\alpha,l)$,
see Marsden et al \cite{mawazh10} and Wang \cite{wa13d}. If we
assume that
$$\mbox{vlift}(u_{(\mu,a)})X_{h_{(\mu,a)}}(\Pi,\Gamma,\alpha,l)
 =(\mathcal{U}_\Pi, \mathcal{U}_\Gamma, \mathcal{U}_\alpha, \mathcal{U}_l)
 \in T_x(\mathcal{O}_{(\mu,a)} \times\mathbb{R}\times\mathbb{R}),
$$ where $x= (\Pi,\Gamma,\alpha,l) \in \mathcal{O}_{(\mu,a)} \times\mathbb{R}\times\mathbb{R}, $
and $\mathcal{U}_\Pi, \; \mathcal{U}_\Gamma \in \mathbb{R}^3, $ and
$\mathcal{U}_\alpha, \; \mathcal{U}_l \in \mathbb{R}. $ Thus, in the
case of non-coincident centers of buoyancy and gravity, the
equations of motion for reduced spacecraft-rotor system with the
control torque $u$ acting on the rotor are given by
\begin{equation}
\left\{\begin{aligned} & \frac{\mathrm{d}\Pi}{\mathrm{d}t}
=\Pi\times\Omega + mgh\Gamma\times\chi + \mathcal{U}_\Pi, \\
& \frac{\mathrm{d}\Gamma}{\mathrm{d}t}=\Gamma\times\Omega + \mathcal{U}_\Gamma, \\
& \frac{\mathrm{d}\alpha}{\mathrm{d}t} =
-\frac{(\Pi_3-l)}{\bar{I}_3}+\frac{l}{J_3} + \mathcal{U}_\alpha, \\
& \frac{\mathrm{d}l}{\mathrm{d}t}= \mathcal{U}_l,
\end{aligned}\right.\label{3.8}\end{equation}

To sum up the above discussion, we have the following theorem.
\begin{theo}
In the case of non-coincident centers of buoyancy and gravity, the
spacecraft-rotor system with the control torque $u$ acting on the
rotor, that is, the 5-tuple $(T^\ast Q,\textmd{SE}(3),\omega_Q,H,u
), $ where $Q= \textmd{SE}(3)\times S^1, $ is a regular point
reducible RCH system. For a point $(\mu,a) \in
\mathfrak{se}^\ast(3)$, the regular value of the momentum map
$\mathbf{J}_Q: \textmd{SE}(3)\times \mathfrak{se}^\ast(3) \times
\mathbb{R} \times \mathbb{R} \to \mathfrak{se}^\ast(3)$, the regular
point reduced system is the 4-tuple $(\mathcal{O}_{(\mu,a)} \times
\mathbb{R} \times \mathbb{R},\tilde{\omega}_{\mathcal{O}_{(\mu,a)}
\times \mathbb{R} \times \mathbb{R}}^{-},h_{(\mu,a)},u_{(\mu,a)}), $
where $\mathcal{O}_{(\mu,a)} \subset \mathfrak{se}^\ast(3)$ is the
coadjoint orbit, $\tilde{\omega}_{\mathcal{O}_{(\mu,a)} \times
\mathbb{R} \times \mathbb{R}}^{-}$ is orbit symplectic form on
$\mathcal{O}_{(\mu,a)} \times \mathbb{R}\times \mathbb{R} $,
$h_{(\mu,a)}(\Pi,\Gamma,\alpha,l)=\pi_{(\mu,a)}(H(A,c,\Pi,\Gamma,\alpha,l))=H(A,c,\Pi,\Gamma,\alpha,l)|_{\mathcal{O}_{(\mu,a)}
\times\mathbb{R}\times\mathbb{R}}$, and
$u_{(\mu,a)}(\Pi,\Gamma,\alpha,l)= \pi_{(\mu,a)}(u(A,c,\Pi,
\\ \Gamma,\alpha,l))= u(A,c,\Pi,
\Gamma,\alpha,l)|_{\mathcal{O}_{(\mu,a)}
\times\mathbb{R}\times\mathbb{R}}$, and its equations of motion are
given by (\ref{3.8}).
\end{theo}

\begin{rema}
When the rigid spacecraft does not carry any internal rotor, in this
case the configuration space is $Q=G=\textmd{SE}(3), $ the motion of
rigid spacecraft is just the rotation motion with drift of a rigid
body, the above symmetric reduction of spacecraft-rotor system is
just the Marsden-Weinstein reduction of a heavy top at a regular
value of momentum map, and the motion equation $(3.8)$ of reduced
spacecraft-rotor system becomes the motion equation of a reduced
heavy top on a coadjoint orbit of Lie group $\textmd{SE}(3)$. See
Marsden and Ratiu \cite{mara99}.
\end{rema}

\section{Hamilton-Jacobi Equation of the Rigid Spacecraft with a Rotor }

It is well-known that Hamilton-Jacobi theory provides a
characterization of the generating functions of certain
time-dependent canonical transformations, such that a given
Hamiltonian system in such a form that its solutions are extremely
easy to find by reduction to the equilibrium, see Abraham and
Marsden \cite{abma78}, Arnold \cite{ar89} and Marsden and Ratiu
\cite{mara99}. In general, we know that it is not easy to find the
solutions of Hamilton's equation. But, if we can get a solution of
Hamilton-Jacobi equation of the Hamiltonian system, by using the
relationship between Hamilton's equation and Hamilton-Jacobi
equation, it is easy to give a special solution of Hamilton's
equation. Thus, it is very important to give explicitly the
Hamilton-Jacobi equation of a Hamiltonian system. Recently, the
author in \cite{wa13d} proved the following Hamilton-Jacobi theorem
for regular point reducible RCH system on the generalization of a
Lie group.

\begin{theo}
For the regular point reducible RCH system $(T^\ast
Q,G,\omega_Q,H,F,W)$ on the generalization of a Lie group $Q=G\times
V$, where $G$ is a Lie group and $V$ is a $k$-dimensional vector
space, assume that $\gamma: Q \rightarrow T^* Q$ is an one-form on
$Q$, and $\gamma^*: T^*T^* Q \rightarrow T^* Q $ is symplectic,
where there is an induced symplectic form $\pi_Q^* \omega_Q$ on
$T^*T^* Q$, $\pi_Q^*: T^* Q \rightarrow T^*T^* Q, \; \pi_Q: T^* Q
\rightarrow Q$, and $\gamma$ is closed with respect to $T\pi_{Q}:
TT^* Q \rightarrow TQ, $ and $\tilde{X}^\gamma = T\pi_{Q}\cdot
\tilde{X} \cdot \gamma$, where $\tilde{X}=X_{(T^\ast
Q,G,\omega_Q,H,F,u)}$ is the dynamical vector field of the regular
point reducible RCH system $(T^*Q,G,\omega_Q,H,F,W)$ with a control
law $u \in W$. Moreover, assume that $\mu \in \mathfrak{g}^\ast $ is
the regular reducible point of the RCH system, and
$\textmd{Im}(\gamma)\subset \mathbf{J}_Q^{-1}(\mu), $ and it is
$G_\mu$-invariant, and $\bar{\gamma}=\pi_\mu(\gamma): Q \rightarrow
\mathcal{O}_\mu \times V\times V^\ast, $ where $G_\mu$ is the
isotropy subgroup of coadjoint action at $\mu$, and $\pi_\mu:
\mathbf{J}_Q^{-1}(\mu) \rightarrow (T^* Q)_\mu \cong \mathcal{O}_\mu
\times V\times V^\ast . $ Then the following two assertions are
equivalent:
 $(\mathrm{i})$ $\tilde{X}^\gamma$ and $\tilde{X}_\mu$ are
$\bar{\gamma}$-related, where $\tilde{X}_\mu= X_{(\mathcal{O}_\mu
\times V\times V^\ast, \tilde{\omega}_{\mathcal{O}_\mu \times V
\times V^\ast }^{-}, h_\mu, f_\mu, u_\mu)}$ is the dynamical vector
field of regular point reduced RCH system $(\mathcal{O}_\mu \times
V\times V^\ast, \tilde{\omega}_{\mathcal{O}_\mu \times V \times
V^\ast }^{-},h_\mu,f_\mu,u_\mu)$; $(\mathrm{ii})$ $X_{h_\mu\cdot
\bar{\gamma}} +\textnormal{vlift}(f_\mu\cdot \bar{\gamma})
+\textnormal{vlift}(u_\mu\cdot \bar{\gamma})=0, $ or $X_{h_\mu\cdot
\bar{\gamma}} +\textnormal{vlift}(f_\mu\cdot \bar{\gamma})
+\textnormal{vlift}(u_\mu\cdot \bar{\gamma})=\tilde{X}^\gamma. $
Moreover, $\gamma$ is a solution of the Hamilton-Jacobi equation $
X_{H\cdot\gamma} +\textnormal{vlift}(F\cdot\gamma)
+\textnormal{vlift}(u\cdot\gamma)=0$, if and only if $\bar{\gamma} $
is a solution of the Hamilton-Jacobi equation $ X_{h_\mu\cdot
\bar{\gamma}} +\textnormal{vlift}(f_\mu\cdot \bar{\gamma})
+\textnormal{vlift}(u_\mu\cdot \bar{\gamma})=0. $
\end{theo}
The maps involved in the above theorem are shown in Diagram-1.
\begin{center}
\hskip 0cm \xymatrix{ \mathbf{J}^{-1}(\mu) \ar[r]^{i_\mu} & T^* Q
\ar[d]_{\pi_Q} \ar[r]^{\pi_Q^*}
& T^*T^* Q \ar[d]_{\gamma^*} & T^*(T^*Q)_\mu \ar[l]_{\pi_\mu^*} \ar[dl]_{\bar{\gamma}^*} \\
  & Q \ar[d]_{\tilde{X}^\gamma} \ar[r]^{\gamma} & T^* Q \ar[d]_{\tilde{X}} \ar[r]^{\pi_\mu} & (T^* Q)_\mu \ar[d]_{\tilde{X}_{\mu}} \\
  & TQ \ar[r]^{T\gamma} & T(T^* Q) \ar[r]^{T\pi_\mu} & T(T^* Q)_\mu   }
\end{center}
$$\mbox{Diagram-1}$$
As an application of the theoretical result, in this section we
regard the rigid spacecraft with an internal rotor as a regular
point reducible RCH system on the generalization of rotation group
$\textmd{SO}(3)\times S^1$ and on the generalization of Euclidean
group $\textmd{SE}(3)\times S^1$, respectively, and give the
Hamilton-Jacobi equations of their reduced RCH systems on the
symplectic leaves by calculation in detail, which show the effect on
controls in Hamilton-Jacobi theory. We shall follow the notations
and conventions introduced in Marsden and Ratiu \cite{mara99},
Marsden \cite{ma92}, Marsden et al \cite{mawazh10}, and Wang
\cite{wa13d}.

\subsection{H-J Equation of Spacecraft-Rotor System with Coincident Centers}

In the following we first give the Hamilton-Jacobi equation for
regular point reduced spacecraft-rotor system with coincident
centers of buoyancy and gravity. From the expression $(2.1)$ of the
Hamiltonian, we know that $H(A,\Pi,\alpha,l)$ is invariant under the
left $\textmd{SO}(3)$-action. For the case $\mu\in
\mathfrak{so}^\ast(3)$ is the regular value of $\mathbf{J}_Q$, we
have the reduced Hamiltonian $h_\mu(\Pi,\alpha,l): \mathcal{O}_\mu
\times\mathbb{R}\times\mathbb{R} \to \mathbb{R}, $ which is given by
$h_\mu(\Pi,\alpha,l)= \pi_\mu(H(A,\Pi,\alpha,l))
=H(A,\Pi,\alpha,l)|_{\mathcal{O}_\mu
\times\mathbb{R}\times\mathbb{R}}, $ and we have the reduced
Hamiltonian vector field
$X_{h_\mu}(K_\mu)=\{K_\mu,h_\mu\}_{-}|_{\mathcal{O}_\mu
\times\mathbb{R}\times\mathbb{R} }. $\\

Assume that $\gamma: \textmd{SO}(3)\times S^1 \rightarrow
T^*(\textmd{SO}(3)\times S^1)$ is an one-form on
$\textmd{SO}(3)\times S^1$, and $\gamma^*: T^*T^* (\textmd{SO}(3)
\times S^1) \rightarrow T^* (\textmd{SO}(3) \times S^1) $ is
symplectic, and $\gamma$ is closed with respect to
$T\pi_{\textmd{SO}(3)\times S^1}: TT^* (\textmd{SO}(3)\times S^1)
\rightarrow T(\textmd{SO}(3)\times S^1), $ and
$\textmd{Im}(\gamma)\subset \mathbf{J}_Q^{-1}(\mu), $ and it is
$\textmd{SO}(3)_\mu$-invariant, and $\bar{\gamma}=\pi_\mu(\gamma):
\textmd{SO}(3)\times S^1 \rightarrow \mathcal{O}_\mu \times
\mathbb{R} \times \mathbb{R}$. Denote by $\bar{\gamma}(A, \alpha)=
(\bar{\gamma}_1,\bar{\gamma}_2,\bar{\gamma}_3,
\bar{\gamma}_4,\bar{\gamma}_5 )(A, \alpha) \in \mathcal{O}_\mu
\times \mathbb{R} \times \mathbb{R}(\subset \mathfrak{so}^\ast(3)
\times \mathbb{R} \times \mathbb{R}), $ that is,
$\Pi=(\Pi_1,\Pi_2,\Pi_3)=(\bar{\gamma}_1,\bar{\gamma}_2,\bar{\gamma}_3)$,
$\alpha=\bar{\gamma}_4, $ and $l=\bar{\gamma}_5 $. Then $h_{\mu}
\cdot \bar{\gamma}: \textmd{SO}(3)\times S^1 \rightarrow \mathbb{R}
$ is given by
\begin{align} h_{\mu} \cdot \bar{\gamma}(A, \alpha)=
H |_{\mathcal{O}_\mu \times \mathbb{R} \times \mathbb{R}} \cdot
\bar{\gamma}(A, \alpha) =
\frac{1}{2}[\frac{\bar{\gamma}_1^2}{\bar{I}_1}+\frac{\bar{\gamma}_2^2}{\bar{I}_2}
+\frac{(\bar{\gamma}_3-
\bar{\gamma}_5)^2}{\bar{I}_3}+\frac{\bar{\gamma}_5^2}{J_3}],
\end{align}
and the vector field
\begin{align*}
& T\bar{\gamma}X_{h_{\mu} \cdot \bar{\gamma}}(\Pi)= X_{h_{\mu}}(\Pi)
\cdot \bar{\gamma} =\{\Pi,h_{\mu}\}_{-}|_{\mathcal{O}_\mu \times
\mathbb{R}\times \mathbb{R}} \cdot \bar{\gamma}(A, \alpha)\\
& =-\Pi\cdot(\nabla_\Pi\Pi\times\nabla_\Pi (h_{\mu})) \cdot
\bar{\gamma}+ \{\Pi,h_{\mu}\}_{\mathbb{R}}|_{\mathcal{O}_\mu \times
\mathbb{R} \times \mathbb{R}} \cdot \bar{\gamma}\\
& =-\nabla_\Pi\Pi\cdot(\nabla_\Pi (h_{\mu})\times \Pi) \cdot
\bar{\gamma}+ (\frac{\partial \Pi}{\partial \alpha} \frac{\partial
h_\mu}{\partial l}- \frac{\partial h_\mu}{\partial
\alpha}\frac{\partial \Pi}{\partial l}) \cdot \bar{\gamma}\\
& =(\Pi_1,\Pi_2,\Pi_3)\times (\frac{\Pi_1}{ \bar{I}_1},\;\;
\frac{\Pi_2}{ \bar{I}_2}, \;\; \frac{(\Pi_3-
l)}{\bar{I}_3}) \cdot \bar{\gamma}\\
&= (\frac{(\bar{I}_2-\bar{I}_3)\bar{\gamma}_2\bar{\gamma}_3-
\bar{I}_2\bar{\gamma}_2\bar{\gamma}_5}{\bar{I}_2\bar{I}_3}, \;
\frac{(\bar{I}_3-\bar{I}_1)\bar{\gamma}_3\bar{\gamma}_1+
\bar{I}_1\bar{\gamma}_1 \bar{\gamma}_5}{\bar{I}_3\bar{I}_1}, \;
\frac{(\bar{I}_1-
\bar{I}_2)\bar{\gamma}_1\bar{\gamma}_2}{\bar{I}_1\bar{I}_2}),
\end{align*}
since $\nabla_\Pi\Pi=1,$ and $\nabla_{\Pi_j} (h_{\mu})= \Pi_j
/\bar{I}_j , \; j= 1,2, \; \nabla_{\Pi_3} (h_{\mu})= (\Pi_3-l)
/\bar{I}_3 , $ and $\frac{\partial \Pi}{\partial \alpha}=
\frac{\partial h_\mu}{\partial \alpha}=0. $

\begin{align*}
& T\bar{\gamma}X_{h_{\mu} \cdot \bar{\gamma}}(\alpha)=
X_{h_{\mu}}(\alpha) \cdot \bar{\gamma} =\{\alpha,
h_{\mu}\}_{-}|_{\mathcal{O}_\mu \times \mathbb{R}\times \mathbb{R}}
\cdot \bar{\gamma}(A, \alpha)\\
& =-\Pi\cdot(\nabla_\Pi\alpha \times\nabla_\Pi (h_{\mu})) \cdot
\bar{\gamma}+ \{\alpha, h_{\mu}\}_{\mathbb{R}}|_{\mathcal{O}_\mu
\times \mathbb{R} \times \mathbb{R}} \cdot \bar{\gamma}\\
& =-\nabla_\Pi\alpha \cdot(\nabla_\Pi (h_{\mu})\times \Pi)\cdot
\bar{\gamma} + (\frac{\partial \alpha}{\partial \alpha}
\frac{\partial h_\mu}{\partial l}- \frac{\partial h_\mu}{\partial
\alpha}\frac{\partial \alpha}{\partial l})\cdot \bar{\gamma} \\
& = -\frac{(\bar{\gamma}_3- \bar{\gamma}_5)}{\bar{I}_3}
+\frac{\bar{\gamma}_5}{J_3},
\end{align*}
since $\nabla_\Pi\alpha=0, \; \frac{\partial \alpha}{\partial l}=0,
$ and $\frac{\partial h_\mu}{\partial l}= -\frac{(\Pi_3-
l)}{\bar{I}_3} +\frac{l}{J_3}$.

\begin{align*}
& T\bar{\gamma}X_{h_{\mu} \cdot \bar{\gamma}}(l)= X_{h_{\mu}}(l)
\cdot \bar{\gamma} =\{l, h_{\mu}\}_{-}|_{\mathcal{O}_\mu \times
\mathbb{R}\times \mathbb{R}} \cdot \bar{\gamma}(A, \alpha)\\
& =-\Pi\cdot(\nabla_\Pi l \times\nabla_\Pi (h_{\mu})) \cdot
\bar{\gamma}+ \{l, h_{\mu}\}_{\mathbb{R}}|_{\mathcal{O}_\mu \times
\mathbb{R} \times \mathbb{R}} \cdot \bar{\gamma}\\
& =-\nabla_\Pi l \cdot(\nabla_\Pi (h_{\mu})\times \Pi)\cdot
\bar{\gamma} + (\frac{\partial l}{\partial \alpha} \frac{\partial
h_\mu}{\partial l}- \frac{\partial h_\mu}{\partial
\alpha}\frac{\partial l}{\partial l})\cdot \bar{\gamma}= 0,
\end{align*}
since $ \nabla_\Pi l=0, $ and $\frac{\partial l}{\partial \alpha}=
\frac{\partial h_\mu}{\partial \alpha}=0. $ If we consider the rigid
spacecraft-rotor system with a control torque $u: T^\ast Q \to W $
acting on the rotor, and $u\in W \subset \mathbf{J}^{-1}_Q(\mu)$ is
invariant under the left $\textmd{SO}(3)$-action, and its reduced
control torque $u_\mu: \mathcal{O}_\mu
\times\mathbb{R}\times\mathbb{R} \to W_\mu $ is given by
$u_\mu(\Pi,\alpha,l)= \pi_\mu(u(A,\Pi,\alpha,l))=
u(A,\Pi,\alpha,l)|_{\mathcal{O}_\mu \times\mathbb{R}\times\mathbb{R}
}, $ where $\pi_\mu: \mathbf{J}_Q^{-1}(\mu) \rightarrow
\mathcal{O}_\mu \times\mathbb{R}\times\mathbb{R}, \; W_\mu=
\pi_\mu(W). $ The dynamical vector field along $\bar{\gamma}$ of
regular point reduced spacecraft-rotor system $(\mathcal{O}_\mu
\times \mathbb{R} \times \mathbb{R},\tilde{\omega}_{\mathcal{O}_\mu
\times \mathbb{R} \times \mathbb{R}}^{-},h_\mu,u_\mu)$ is given by
$$ X_{(\mathcal{O}_\mu \times \mathbb{R} \times
\mathbb{R},\tilde{\omega}_{\mathcal{O}_\mu \times \mathbb{R} \times
\mathbb{R}}^{-},h_\mu,u_\mu)} \cdot \bar{\gamma}= X_{h_\mu} \cdot
\bar{\gamma}+ \mbox{vlift}(u_\mu) \cdot \bar{\gamma}, $$ where
$\mbox{vlift}(u_\mu) \cdot \bar{\gamma} = T\bar{\gamma}\cdot
\mbox{vlift}(u_\mu \cdot \bar{\gamma})= T\bar{\gamma}\cdot
\mbox{vlift}(u_\mu \cdot \bar{\gamma})X_{h_\mu \cdot \bar{\gamma}}
\in T(\mathcal{O}_\mu \times\mathbb{R}\times\mathbb{R}). $ Assume
that
$$\mbox{vlift}(u_\mu) \cdot \bar{\gamma}(A, \alpha)= (U_1, U_2, U_3,
U_4, U_5 )(A, \alpha) \in T_x(\mathcal{O}_\mu
\times\mathbb{R}\times\mathbb{R}), $$ where $x= (\Pi,\alpha,l)\in
\mathcal{O}_\mu \times\mathbb{R}\times\mathbb{R}, $ and
$U_i(A,\alpha): \textmd{SO}(3)\times S^1 \rightarrow \mathbb{R}, \;
i=1,2, \cdots,5, $ then from the Hamilton-Jacobi equations
$X_{h_{\mu}\cdot \bar{\gamma}}+ \mbox{vlift}(u_{\mu}\cdot
\bar{\gamma})=0, $ we have that
\begin{align*}
0 &= T\bar{\gamma}(X_{h_{\mu}\cdot \bar{\gamma}}+
\mbox{vlift}(u_{\mu}\cdot \bar{\gamma}))= X_{h_{\mu}}\cdot
\bar{\gamma}+ \mbox{vlift}(u_{\mu})\cdot \bar{\gamma}\\ &=
X_{(\mathcal{O}_{\mu} \times \mathbb{R} \times
\mathbb{R},\tilde{\omega}_{\mathcal{O}_{\mu} \times \mathbb{R}
\times \mathbb{R}}^{-},h_{\mu},u_{\mu})} \cdot \bar{\gamma},
\end{align*}
that is, the dynamical vector field along $\bar{\gamma}$ of regular
point reduced spacecraft-rotor system $(\mathcal{O}_{\mu} \times
\mathbb{R} \times \mathbb{R},\tilde{\omega}_{\mathcal{O}_{\mu}
\times \mathbb{R} \times \mathbb{R}}^{-},h_{\mu},u_{\mu})$ is its
equilibrium. Thus, in the case of coincident centers of buoyancy and
gravity, the Hamilton-Jacobi equations for reduced spacecraft-rotor
system with the control torque $u$ acting on the rotor are given by
\begin{equation}
 \left\{\begin{aligned}
& \frac{(\bar{I}_2-\bar{I}_3)\bar{\gamma}_2\bar{\gamma}_3-
\bar{I}_2\bar{\gamma}_2\bar{\gamma}_5}{\bar{I}_2\bar{I}_3}+ U_1 =
0,\\
& \frac{(\bar{I}_3-\bar{I}_1)\bar{\gamma}_3\bar{\gamma}_1+
\bar{I}_1\bar{\gamma}_1 \bar{\gamma}_5}{\bar{I}_3\bar{I}_1}+ U_2 =
0,\\
& \frac{(\bar{I}_1-
\bar{I}_2)\bar{\gamma}_1\bar{\gamma}_2}{\bar{I}_1\bar{I}_2}+ U_3 =
0,\\
& -\frac{(\bar{\gamma}_3- \bar{\gamma}_5)}{\bar{I}_3}
+\frac{\bar{\gamma}_5}{J_3} + U_4 = 0,
\\ & U_5 =0.
\end{aligned} \right. \label{4.2}
\end{equation}

To sum up the above discussion, we have the following theorem.
\begin{theo}
In the case of coincident centers of buoyancy and gravity, for a
point $\mu \in \mathfrak{so}^\ast(3)$, the regular value of the
momentum map $\mathbf{J}_Q: \textmd{SO}(3)\times
\mathfrak{so}^\ast(3) \times \mathbb{R} \times \mathbb{R} \to
\mathfrak{so}^\ast(3)$, the regular point reduced system of
spacecraft-rotor system with the control torque $u$ acting on the
rotor $(T^\ast Q,\textmd{SO}(3),\omega_Q,H,u), $ where $Q=
\textmd{SO}(3)\times S^1, $ is the 4-tuple $(\mathcal{O}_\mu \times
\mathbb{R} \times \mathbb{R},\tilde{\omega}_{\mathcal{O}_\mu \times
\mathbb{R} \times \mathbb{R}}^{-},h_\mu,u_\mu), $ where
$\mathcal{O}_\mu \subset \mathfrak{so}^\ast(3)$ is the coadjoint
orbit, $\tilde{\omega}_{\mathcal{O}_\mu \times \mathbb{R} \times
\mathbb{R}}^{-}$ is orbit symplectic form on $\mathcal{O}_\mu \times
\mathbb{R}\times \mathbb{R} $, $h_\mu(\Pi,\alpha,l)=
\pi_\mu(H(A,\Pi,\alpha,l))=H(A,\Pi,\alpha,l)|_{\mathcal{O}_\mu
\times\mathbb{R}\times\mathbb{R}}$, $u_\mu(\Pi,\alpha,l)=
\pi_\mu(u(A,\Pi,\alpha,l))= u(A,\Pi,\alpha,l)|_{\mathcal{O}_\mu
\times\mathbb{R}\times\mathbb{R}}$. Assume that $\gamma:
\textmd{SO}(3)\times S^1 \rightarrow T^*(\textmd{SO}(3)\times S^1)$
is an one-form on $\textmd{SO}(3)\times S^1$, and $\gamma^*: T^*T^*
(\textmd{SO}(3) \times S^1) \rightarrow T^* (\textmd{SO}(3) \times
S^1) $ is symplectic, and $\gamma$ is closed with respect to
$T\pi_Q: TT^* Q \rightarrow TQ, $ and $\textmd{Im}(\gamma)\subset
\mathbf{J}_Q^{-1}(\mu), $ and it is $\textmd{SO}(3)_\mu$-invariant,
where $\textmd{SO}(3)_\mu$ is the isotropy subgroup of coadjoint
$\textmd{SO}(3)$-action at the point $\mu\in\mathfrak{so}^\ast(3)$,
and $\bar{\gamma}=\pi_\mu(\gamma): \textmd{SO}(3)\times S^1
\rightarrow \mathcal{O}_\mu \times \mathbb{R}\times \mathbb{R}$.
Then $\bar{\gamma} $ is a solution of either Hamilton-Jacobi
equation of reduced spacecraft-rotor system given by $(4.2)$, or the
equation $X_{h_\mu\cdot \bar{\gamma}} +\textnormal{vlift}(u_\mu\cdot
\bar{\gamma})=\tilde{X}^\gamma,$ if and only if $\tilde{X}^\gamma$
and $\tilde{X}_\mu$ are $\bar{\gamma}$-related, where
$\tilde{X}^\gamma= T\pi_{Q} \cdot \tilde{X} \cdot \gamma, \;
\tilde{X}= X_{(T^\ast Q,\textmd{SO}(3),\omega_Q,H,u)}, $ and
$\tilde{X}_\mu= X_{(\mathcal{O}_\mu \times \mathbb{R} \times
\mathbb{R}, \tilde{\omega}_{\mathcal{O}_\mu \times \mathbb{R} \times
\mathbb{R} }^{-}, h_\mu, u_\mu)}$.
\end{theo}

\begin{rema}
When the rigid spacecraft does not carry any internal rotor, in this
case the configuration space is $Q=G=\textmd{SO}(3), $ the above
Hamilton-Jacobi equation $(4.2)$ of reduced spacecraft-rotor system
is just the Lie-Poisson Hamilton-Jacobi equation of
Marsden-Weinstein reduced Hamiltonian system $(\mathcal{O}_\mu,
\omega_{\mathcal{O}_\mu}^{-}, h_{\mathcal{O}_\mu} )$ on the
coadjoint orbit of Lie group $\textmd{SO}(3)$, since the symplecitic
structure on the coadjoint orbit $\mathcal{O}_\mu$ is induced by the
(-)-Lie-Poisson brackets on $\mathfrak{so}^\ast(3)$. See Marsden and
Ratiu \cite{mara99}, Ge and Marsden \cite{gema88}, and Wang
\cite{wa13a}.
\end{rema}

\subsection{H-J Equation of Spacecraft-Rotor System with Non-coincident Centers}

In the following we shall give the Hamilton-Jacobi equation for
regular point reduced spacecraft-rotor system with non-coincident
centers of buoyancy and gravity. From the expression $(2.2)$ of the
Hamiltonian, we know that $H(A,c,\Pi,\Gamma,\alpha,l)$ is invariant
under the left $\textmd{SE}(3)$-action. For the case $(\mu,a)\in
\mathfrak{se}^\ast(3)$ is the regular value of $\mathbf{J}_Q$, we
have the reduced Hamiltonian
$h_{(\mu,a)}(\Pi,\Gamma,\alpha,l):\mathcal{O}_{(\mu,a)}\times\mathbb{R}
\times \mathbb{R} (\subset \mathfrak{se}^\ast (3)\times
\mathbb{R}\times \mathbb{R}) \to \mathbb{R}, $ which is given by
$h_{(\mu,a)}(\Pi,\Gamma,\alpha,l)
=\pi_{(\mu,a)}(H(A,c,\Pi,\Gamma,\alpha,l))=H(A,c,\Pi,\Gamma,\alpha,l)|_{\mathcal{O}_{(\mu,a)}\times
\mathbb{R}\times \mathbb{R}}$, and we have the reduced Hamiltonian
vector field
$X_{h_{(\mu,a)}}(K_{(\mu,a)})=\{K_{(\mu,a)},h_{(\mu,a)}\}_{-}|_{\mathcal{O}_{(\mu,a)}
\times \mathbb{R}\times \mathbb{R}}. $\\

Assume that $\gamma: \textmd{SE}(3)\times S^1 \rightarrow T^*
(\textmd{SE}(3) \times S^1) $ is an one-form on $\textmd{SE}(3)
\times S^1 $, and $\gamma^*: T^*T^* (\textmd{SE}(3) \times S^1)
\rightarrow T^* (\textmd{SE}(3) \times S^1) $ is symplectic, and
$\gamma$ is closed with respect to $T\pi_{\textmd{SE}(3)\times S^1}:
TT^* (\textmd{SE}(3) \times S^1) \rightarrow T(\textmd{SE}(3) \times
S^1), $ and $\textmd{Im}(\gamma)\subset \mathbf{J}^{-1}((\mu,a)), $
and it is $\textmd{SE}(3)_{(\mu,a)}$-invariant, and
$\bar{\gamma}=\pi_{(\mu,a)}(\gamma): \textmd{SE}(3) \times S^1
\rightarrow \mathcal{O}_{(\mu,a)}\times \mathbb{R} \times
\mathbb{R}$. Denote by $\bar{\gamma}(A,c,\alpha)= (\bar{\gamma}_1,
\bar{\gamma}_2, \bar{\gamma}_3, \bar{\gamma}_4, \bar{\gamma}_5,
\bar{\gamma}_6, \bar{\gamma}_7,\\ \bar{\gamma}_8 )(A,c,\alpha) \in
\mathcal{O}_{(\mu,a)} \times \mathbb{R} \times \mathbb{R} (\subset
\mathfrak{se}^\ast(3) \times \mathbb{R} \times \mathbb{R}), $ that
is,
$\Pi=(\Pi_1,\Pi_2,\Pi_3)=(\bar{\gamma}_1,\bar{\gamma}_2,\bar{\gamma}_3)$,
$\Gamma=(\Gamma_1,\Gamma_2,\Gamma_3)=(\bar{\gamma}_4,\bar{\gamma}_5,\bar{\gamma}_6)$,
$\alpha=\bar{\gamma}_7, $ and $l=\bar{\gamma}_8 $. Then $h_{(\mu,a)}
\cdot \bar{\gamma}: \textmd{SE}(3) \times S^1 \rightarrow \mathbb{R}
$ is given by
\begin{align}
& h_{(\mu,a)} \cdot \bar{\gamma}(A,c,\alpha) = H
|_{\mathcal{O}_{(\mu,a)} \times \mathbb{R} \times \mathbb{R}}\cdot
\bar{\gamma}(A,c,\alpha) \nonumber \\ &
=\frac{1}{2}[\frac{\bar{\gamma}_1^2}{\bar{I}_1}+\frac{\bar{\gamma}_2^2}{\bar{I}_2}
+\frac{(\bar{\gamma}_3-\bar{\gamma}_8)^2}{\bar{I}_3}+\frac{\bar{\gamma}_8^2}{J_3}]
+ mgh(\bar{\gamma}_4\cdot\chi_1+ \bar{\gamma}_5\cdot\chi_2+
\bar{\gamma}_6\cdot\chi_3),
\end{align} and the vector field
\begin{align*}
& T\bar{\gamma}X_{h_{(\mu,a)} \cdot \bar{\gamma}}(\Pi)=
X_{h_{(\mu,a)}}(\Pi) \cdot \bar{\gamma}= \{\Pi,\;
h_{(\mu,a)}\}_{-}|_{\mathcal{O}_{(\mu,a)}
\times \mathbb{R}\times \mathbb{R}}\cdot \bar{\gamma}(A,c,\alpha)\\
& =-\Pi\cdot(\nabla_\Pi\Pi\times\nabla_\Pi (h_{(\mu,a)})) \cdot
\bar{\gamma}-\Gamma\cdot(\nabla_\Pi\Pi\times\nabla_\Gamma
(h_{(\mu,a)})-\nabla_\Pi
(h_{(\mu,a)}) \times\nabla_\Gamma\Pi)\cdot \bar{\gamma}\\
& \;\;\;\; + \{\Pi,\;
h_{(\mu,a)}\}_{\mathbb{R}}|_{\mathcal{O}_{(\mu,a)} \times
\mathbb{R}\times \mathbb{R}}\cdot \bar{\gamma}\\
& = -\nabla_\Pi\Pi \cdot(\nabla_\Pi (h_{(\mu,a)}) \times \Pi)\cdot
\bar{\gamma}- \nabla_\Pi\Pi \cdot (\nabla_\Gamma
(h_{(\mu,a)}) \times \Gamma)\cdot \bar{\gamma}\\
& \;\;\;\; + (\frac{\partial \Pi}{\partial \alpha} \frac{\partial
h_{(\mu,a)}}{\partial l}- \frac{\partial h_{(\mu,a)}}{\partial
\alpha}\frac{\partial \Pi}{\partial l})\cdot \bar{\gamma}\\
&=(\Pi_1,\Pi_2,\Pi_3)\times (\frac{\Pi_1}{ \bar{I}_1},
\frac{\Pi_2}{\bar{I}_2}, \frac{(\Pi_3- l)}{\bar{I}_3})\cdot
\bar{\gamma}
+mgh(\Gamma_1,\Gamma_2,\Gamma_3)\times (\chi_1,\chi_2,\chi_3)\cdot \bar{\gamma}\\
&= ( \; \frac{(\bar{I}_2-\bar{I}_3)\bar{\gamma}_2\bar{\gamma}_3-
\bar{I}_2\bar{\gamma}_2 \bar{\gamma}_8}{\bar{I}_2\bar{I}_3}+
mgh(\bar{\gamma}_5\chi_3-\bar{\gamma}_6\chi_2), \\ & \;\;\;\;\;\;
\frac{(\bar{I}_3-\bar{I}_1)\bar{\gamma}_3 \bar{\gamma}_1+
\bar{I}_1\bar{\gamma}_1\bar{\gamma}_8 }{\bar{I}_3\bar{I}_1} +
mgh(\bar{\gamma}_6\chi_1-\bar{\gamma}_4\chi_3), \;\;\;\;\;\;
\frac{(\bar{I}_1-\bar{I}_2)\bar{\gamma}_1
\bar{\gamma}_2}{\bar{I}_1\bar{I}_2} +
mgh(\bar{\gamma}_4\chi_2-\bar{\gamma}_5\chi_1) \; ),
\end{align*}
since $\nabla_\Pi\Pi=1, \; \nabla_\Gamma\Pi =0, \; \frac{\partial
\Pi}{\partial \alpha}= \frac{\partial \Pi}{\partial l}=0, \;
\Gamma=(\Gamma_1,\Gamma_2,\Gamma_3), \; \chi=(\chi_1,\chi_2,\chi_3),
$ and $\nabla_{\Pi_j} (h_{(\mu,a)})= \Pi_j /\bar{I}_j , \; j= 1,2, $
$ \nabla_{\Pi_3} (h_{(\mu,a)})= (\Pi_3- l)/\bar{I}_3, \;
\nabla_{\Gamma_j} h_{(\mu,a)}= mgh\chi_j, \; j=1,2,3. $

\begin{align*}
& T\bar{\gamma}X_{h_{(\mu,a)} \cdot \bar{\gamma}}(\Gamma)=
X_{h_{(\mu,a)}}(\Gamma) \cdot \bar{\gamma}= \{\Gamma,\;
h_{(\mu,a)}\}_{-}|_{\mathcal{O}_{(\mu,a)} \times \mathbb{R}\times
\mathbb{R}}\cdot \bar{\gamma}(A,c,\alpha)\\
& =-\Pi\cdot(\nabla_\Pi\Gamma\times\nabla_\Pi (h_{(\mu,a)})) \cdot
\bar{\gamma}-\Gamma\cdot(\nabla_\Pi\Gamma\times\nabla_\Gamma
(h_{(\mu,a)})-\nabla_\Pi
(h_{(\mu,a)}) \times\nabla_\Gamma\Gamma)\cdot \bar{\gamma}\\
& \;\;\;\; + \{\Gamma,\;
h_{(\mu,a)}\}_{\mathbb{R}}|_{\mathcal{O}_{(\mu,a)} \times
\mathbb{R}\times \mathbb{R}}\cdot \bar{\gamma}\\
& =\nabla_\Gamma\Gamma\cdot(\Gamma\times\nabla_\Pi
(h_{(\mu,a)}))\cdot \bar{\gamma}+ (\frac{\partial \Gamma}{\partial
\alpha} \frac{\partial h_{(\mu,a)}}{\partial l}- \frac{\partial
h_{(\mu,a)}}{\partial
\alpha}\frac{\partial \Gamma}{\partial l})\cdot \bar{\gamma}\\
&=(\Gamma_1,\Gamma_2,\Gamma_3)\times (\frac{\Pi_1}{ \bar{I}_1},
\frac{\Pi_2}{ \bar{I}_2}, \frac{(\Pi_3- l)}{\bar{I}_3})\cdot
\bar{\gamma} \\
&= ( \; \frac{\bar{I}_2\bar{\gamma}_5(\bar{\gamma}_3-
\bar{\gamma}_8)- \bar{I}_3\bar{\gamma}_6
\bar{\gamma}_2}{\bar{I}_2\bar{I}_3}, \;\;\;
\frac{\bar{I}_3\bar{\gamma}_6 \bar{\gamma}_1-
\bar{I}_1\bar{\gamma}_4(\bar{\gamma}_3-
\bar{\gamma}_8)}{\bar{I}_3\bar{I}_1}, \;\;\;
\frac{\bar{I}_1\bar{\gamma}_4 \bar{\gamma}_2-
\bar{I}_2\bar{\gamma}_5 \bar{\gamma}_1}{\bar{I}_1\bar{I}_2} \; ),
\end{align*}
since $\nabla_\Gamma \Gamma =1, \; \nabla_\Pi\Gamma =0, \;
\frac{\partial \Gamma}{\partial \alpha}= \frac{\partial
\Gamma}{\partial l}=0. $

\begin{align*}
& T\bar{\gamma}X_{h_{(\mu,a)} \cdot \bar{\gamma}}(\alpha)=
X_{h_{(\mu,a)}}(\alpha) \cdot \bar{\gamma}= \{\alpha,\;
h_{(\mu,a)}\}_{-}|_{\mathcal{O}_{(\mu,a)} \times \mathbb{R}\times
\mathbb{R}}\cdot \bar{\gamma}(A,c,\alpha)\\
& =-\Pi\cdot(\nabla_\Pi \alpha \times\nabla_\Pi (h_{(\mu,a)})) \cdot
\bar{\gamma}-\Gamma\cdot(\nabla_\Pi \alpha \times\nabla_\Gamma
(h_{(\mu,a)})-\nabla_\Pi
(h_{(\mu,a)}) \times\nabla_\Gamma \alpha)\cdot \bar{\gamma}\\
& \;\;\;\; + \{\alpha,\;
h_{(\mu,a)}\}_{\mathbb{R}}|_{\mathcal{O}_{(\mu,a)} \times
\mathbb{R}\times \mathbb{R}}\cdot \bar{\gamma}\\
& = (\frac{\partial \alpha}{\partial \alpha} \frac{\partial
h_{(\mu,a)}}{\partial l}- \frac{\partial h_{(\mu,a)}}{\partial
\alpha}\frac{\partial \alpha}{\partial l})\cdot \bar{\gamma} =
-\frac{(\bar{\gamma}_3- \bar{\gamma}_8)}{ \bar{I}_3}+
\frac{\bar{\gamma}_8}{J_3},
\end{align*}
since $\nabla_\Pi\alpha= \nabla_\Gamma \alpha=0, \; \frac{\partial
\alpha}{\partial l}=0, $ and $\frac{\partial h_{(\mu,a)}}{\partial
l}= -\frac{(\Pi_3- l)}{\bar{I}_3} +\frac{l}{J_3}$.

\begin{align*}
& T\bar{\gamma}X_{h_{(\mu,a)} \cdot \bar{\gamma}}(l)=
X_{h_{(\mu,a)}}(l) \cdot \bar{\gamma}= \{l,\;
h_{(\mu,a)}\}_{-}|_{\mathcal{O}_{(\mu,a)} \times \mathbb{R}\times
\mathbb{R}}\cdot \bar{\gamma}(A,c,\alpha)\\
& =-\Pi\cdot(\nabla_\Pi l \times\nabla_\Pi (h_{(\mu,a)})) \cdot
\bar{\gamma}-\Gamma\cdot(\nabla_\Pi l \times\nabla_\Gamma
(h_{(\mu,a)})-\nabla_\Pi
(h_{(\mu,a)}) \times\nabla_\Gamma l)\cdot \bar{\gamma}\\
& \;\;\;\; + \{l,\;
h_{(\mu,a)}\}_{\mathbb{R}}|_{\mathcal{O}_{(\mu,a)} \times
\mathbb{R}\times \mathbb{R}}\cdot \bar{\gamma}\\
& = (\frac{\partial l}{\partial \alpha} \frac{\partial
h_{(\mu,a)}}{\partial l}- \frac{\partial h_{(\mu,a)}}{\partial
\alpha}\frac{\partial l}{\partial l})\cdot \bar{\gamma}= 0,
\end{align*}
since $ \nabla_\Pi l= \nabla_\Gamma l =0, \; \frac{\partial
l}{\partial \alpha}= 0, $ and $\frac{\partial h_{(\mu,a)} }{\partial
\alpha}= 0. $ If we consider the spacecraft-rotor system with a
control torque $u: T^\ast Q \to W $ acting on the rotor, and $u\in W
\subset \mathbf{J}^{-1}_Q((\mu,a))$ is invariant under the left
$\textmd{SE}(3)$-action, and its reduced control torque
$u_{(\mu,a)}: \mathcal{O}_{(\mu,a)} \times\mathbb{R}\times\mathbb{R}
\to W_{(\mu,a)} $ is given by
$u_{(\mu,a)}(\Pi,\Gamma,\alpha,l)=\pi_{(\mu,a)}(u(A,c,\Pi,\Gamma,\alpha,l))=
u(A,c,\Pi,\Gamma,\alpha,l)|_{\mathcal{O}_{(\mu,a)}
\times\mathbb{R}\times\mathbb{R} }, $ where $\pi_{(\mu,a)}:
\mathbf{J}_Q^{-1}((\mu,a)) \rightarrow \mathcal{O}_{(\mu,a)}
\times\mathbb{R}\times\mathbb{R}, \, ; W_{(\mu,a)}=
\pi_{(\mu,a)}(W). $ The dynamical vector field along $\bar{\gamma}$
of regular point reduced spacecraft-rotor system
$(\mathcal{O}_{(\mu,a)} \times \mathbb{R} \times
\mathbb{R},\tilde{\omega}_{\mathcal{O}_{(\mu,a)} \times \mathbb{R}
\times \mathbb{R}}^{-},h_{(\mu,a)},u_{(\mu,a)})$ is given by
$$ X_{(\mathcal{O}_{(\mu,a)} \times \mathbb{R} \times
\mathbb{R},\tilde{\omega}_{\mathcal{O}_{(\mu,a)} \times \mathbb{R}
\times \mathbb{R}}^{-},h_{(\mu,a)},u_{(\mu,a)})} \cdot \bar{\gamma}=
X_{h_{(\mu,a)}} \cdot \bar{\gamma}+ \mbox{vlift}(u_{(\mu,a)}) \cdot
\bar{\gamma}, $$ where $\mbox{vlift}(u_{(\mu,a)}) \cdot \bar{\gamma}
= T\bar{\gamma} \cdot \mbox{vlift}(u_{(\mu,a)}\cdot \bar{\gamma}) =
T\bar{\gamma} \cdot \mbox{vlift}(u_{(\mu,a)}\cdot
\bar{\gamma})X_{h_{(\mu,a)}\cdot \bar{\gamma}} \in
T(\mathcal{O}_{(\mu,a)} \times\mathbb{R}\times\mathbb{R} ). $ Assume
that
$$\mbox{vlift}(u_{(\mu,a)}) \cdot \bar{\gamma}(A, c,\alpha)= (U_1,
U_2, U_3, U_4, U_5, U_6, U_7, U_8 )(A,c,\alpha) \in
T_x(\mathcal{O}_{(\mu,a)} \times\mathbb{R}\times\mathbb{R} ), $$
where $x= (\Pi,\Gamma,\alpha,l) \in \mathcal{O}_{(\mu,a)}
\times\mathbb{R}\times\mathbb{R}, $ and $U_i(A,c,\alpha):
\textmd{SE}(3)\times S^1 \rightarrow \mathbb{R}, \; i= 1,2,\cdots,8,
$ then from the Hamilton-Jacobi equations $X_{h_{(\mu,a)}\cdot
\bar{\gamma}}+ \mbox{vlift}(u_{(\mu,a)}\cdot \bar{\gamma})=0, $ we
have that
\begin{align*}
0 &= T\bar{\gamma}(X_{h_{(\mu,a)}\cdot \bar{\gamma}}+
\mbox{vlift}(u_{(\mu,a)}\cdot \bar{\gamma}))= X_{h_{(\mu,a)}}\cdot
\bar{\gamma}+ \mbox{vlift}(u_{(\mu,a)})\cdot \bar{\gamma}
\\ &= X_{(\mathcal{O}_{(\mu,a)} \times \mathbb{R} \times
\mathbb{R},\tilde{\omega}_{\mathcal{O}_{(\mu,a)} \times \mathbb{R}
\times \mathbb{R}}^{-},h_{(\mu,a)},u_{(\mu,a)})} \cdot \bar{\gamma},
\end{align*}
that is, the dynamical vector field along $\bar{\gamma}$ of regular
point reduced spacecraft-rotor system $(\mathcal{O}_{(\mu,a)} \times
\mathbb{R} \times \mathbb{R},\tilde{\omega}_{\mathcal{O}_{(\mu,a)}
\times \mathbb{R} \times \mathbb{R}}^{-},h_{(\mu,a)},u_{(\mu,a)})$
is its equilibrium. Thus, in the case of non-coincident centers of
buoyancy and gravity, the Hamilton-Jacobi equations for reduced
spacecraft-rotor system with the control torque $u$ acting on the
rotor are given by
\begin{equation}
\left\{\begin{aligned} &
\frac{(\bar{I}_2-\bar{I}_3)\bar{\gamma}_2\bar{\gamma}_3-
\bar{I}_2\bar{\gamma}_2 \bar{\gamma}_8}{\bar{I}_2\bar{I}_3}+
mgh(\bar{\gamma}_5\chi_3-\bar{\gamma}_6\chi_2)+ U_1=0,\\
& \frac{(\bar{I}_3-\bar{I}_1)\bar{\gamma}_3 \bar{\gamma}_1+
\bar{I}_1\bar{\gamma}_1\bar{\gamma}_8 }{\bar{I}_3\bar{I}_1} +
mgh(\bar{\gamma}_6\chi_1-\bar{\gamma}_4\chi_3)+ U_2=0,\\
& \frac{(\bar{I}_1-\bar{I}_2)\bar{\gamma}_1
\bar{\gamma}_2}{\bar{I}_1\bar{I}_2} +
mgh(\bar{\gamma}_4\chi_2-\bar{\gamma}_5\chi_1)+ U_3 =0,\\
& \frac{\bar{I}_2\bar{\gamma}_5(\bar{\gamma}_3- \bar{\gamma}_8)-
\bar{I}_3\bar{\gamma}_6 \bar{\gamma}_2}{\bar{I}_2\bar{I}_3}+ U_4 =0,\\
& \frac{\bar{I}_3\bar{\gamma}_6 \bar{\gamma}_1-
\bar{I}_1\bar{\gamma}_4(\bar{\gamma}_3-
\bar{\gamma}_8)}{\bar{I}_3\bar{I}_1}+ U_5=0, \\
& \frac{\bar{I}_1\bar{\gamma}_4 \bar{\gamma}_2-
\bar{I}_2\bar{\gamma}_5 \bar{\gamma}_1}{\bar{I}_1\bar{I}_2}+ U_6 =0,\\
& -\frac{(\bar{\gamma}_3- \bar{\gamma}_8)}{ \bar{I}_3}+
\frac{\bar{\gamma}_8}{J_3} + U_7 =0,
\\ & U_8 =0 .
\end{aligned} \right. \label{4.4}
\end{equation}

To sum up the above discussion, we have the following theorem.
\begin{theo}
In the case of non-coincident centers of buoyancy and gravity, for a
point $(\mu,a)\in \mathfrak{se}^\ast(3)$, the regular value of the
momentum map $\mathbf{J}_Q: \textmd{SE}(3)\times
\mathfrak{se}^\ast(3)\times \mathbb{R}\times \mathbb{R} \to
\mathfrak{se}^\ast(3)$, the regular point reduced system of
spacecraft-rotor system with the control torque $u$ acting on the
rotor $(T^\ast Q,\textmd{SE}(3),\omega_Q,H,u), $ where $Q=
\textmd{SE}(3)\times S^1, $ is the 4-tuple
$(\mathcal{O}_{(\mu,a)}\times \mathbb{R}\times \mathbb{R},
\tilde{\omega}^{-}_{\mathcal{O}_{(\mu,a)}\times \mathbb{R}\times
\mathbb{R}},h_{(\mu,a)},u_{(\mu,a)})$, where $\mathcal{O}_{(\mu,a)}
\subset \mathfrak{se}^\ast(3)$ is the coadjoint orbit,
$\tilde{\omega}^{-}_{\mathcal{O}_{(\mu,a)}\times \mathbb{R}\times
\mathbb{R}}$ is orbit symplectic form on $\mathcal{O}_{(\mu,a)}
\times \mathbb{R}\times \mathbb{R} $,
$h_{(\mu,a)}(\Pi,\Gamma,\alpha,l)=\pi_{(\mu,a)}(H(A,c,\Pi,\Gamma,\alpha,l))
=H(A,c,\Pi,\Gamma,\alpha,l)|_{\mathcal{O}_{(\mu,a)}\times
\mathbb{R}\times \mathbb{R}}, $ and
$u_{(\mu,a)}(\Pi,\Gamma,\alpha,l)=\pi_{(\mu,a)}(u(A,c,\Pi,\Gamma,\alpha,l))=
u(A,c,\Pi,\Gamma,\alpha,l)|_{\mathcal{O}_{(\mu,a)}
\times\mathbb{R}\times\mathbb{R} }. $ Assume that $\gamma:
\textmd{SE}(3) \times S^1 \rightarrow T^* (\textmd{SE}(3) \times
S^1) $ is an one-form on $\textmd{SE}(3) \times S^1 $, and
$\gamma^*: T^*T^* (\textmd{SE}(3) \times S^1) \rightarrow T^*
(\textmd{SE}(3) \times S^1) $ is symplectic, and $\gamma$ is closed
with respect to $T\pi_{Q}: TT^* Q \rightarrow TQ, $ and
$\textmd{Im}(\gamma)\subset \mathbf{J}^{-1}(\mu,a), $ and it is
$\textmd{SE}(3)_{(\mu,a)}$-invariant, where
$\textmd{SE}(3)_{(\mu,a)}$ is the isotropy subgroup of coadjoint
$\textmd{SE}(3)$-action at the point
$(\mu,a)\in\mathfrak{se}^\ast(3)$, and
$\bar{\gamma}=\pi_{(\mu,a)}(\gamma): \textmd{SE}(3) \times S^1
\rightarrow \mathcal{O}_{(\mu,a)} \times \mathbb{R} \times
\mathbb{R}$. Then $\bar{\gamma} $ is a solution of either the
Hamilton-Jacobi equation of reduced spacecraft-rotor system given by
$(4.4)$, or the equation
$X_{h_{\mathcal{O}_{(\mu,a)}}\cdot\bar{\gamma}}+
\textnormal{vlift}(u_{(\mu,a)}\cdot\bar{\gamma})= \tilde{X}^\gamma
,$ if and only if $\tilde{X}^\gamma$ and $\tilde{X}_{(\mu,a)}$ are
$\bar{\gamma}$-related, where $\tilde{X}^\gamma = T\pi_{Q}\cdot
\tilde{X} \cdot \gamma, \; \tilde{X}= X_{(T^\ast
Q,\textmd{SE}(3),\omega_Q,H,u)}, $ and $\tilde{X}_{(\mu,a)}=
X_{(\mathcal{O}_{(\mu,a)} \times \mathbb{R}\times \mathbb{R},
\tilde{\omega}_{\mathcal{O}_{(\mu,a)} \times \mathbb{R} \times
\mathbb{R} }^{-}, h_{(\mu,a)}, u_{(\mu,a)})}$.
\end{theo}
\begin{rema}
When the rigid spacecraft does not carry any internal rotor, in this
case the configuration space is $Q=G=\textmd{SE}(3), $ the above
Hamilton-Jacobi equation $(4.4)$ of reduced spacecraft-rotor system
is just the Lie-Poisson Hamilton-Jacobi equation of
Marsden-Weinstein reduced Hamiltonian system
$(\mathcal{O}_{(\mu,a)}, \omega_{\mathcal{O}_{(\mu,a)}}^{-},
h_{\mathcal{O}_{(\mu,a)}} )$ on the coadjoint orbit of Lie group
$\textmd{SE}(3)$, since the symplecitic structure on the coadjoint
orbit $\mathcal{O}_{(\mu,a)}$ is induced by the (-)-Lie-Poisson
brackets on $\mathfrak{se}^\ast(3)$. See Marsden and Ratiu
\cite{mara99}, Ge and Marsden \cite{gema88}, and Wang \cite{wa13a}.
\end{rema}

The theory of controlled mechanical system is a very important
subject, following the theoretical and applied development of
geometric mechanics, a lot of important problems about this subject
are being explored and studied. In this paper, as an application of
the symplectic reduction and Hamilton-Jacobi theory of regular
controlled Hamiltonian systems with symmetry, in the cases of
coincident and non-coincident centers of buoyancy and gravity, we
give explicitly the motion equation and Hamilton-Jacobi equation of
reduced spacecraft-rotor system on a symplectic leaf by calculation
in detail, respectively, which show the effect on controls in
regular symplectic reduction and Hamilton-Jacobi theory. We also
note that in \cite{mlml08}, the authors study the dynamics of a
rigid spacecraft under the influence of gravity torques and give the
dynamical equations in a first-order form with a special coefficient
matrix. Moreover, we hope to study the stabilization of rigid
spacecraft with an internal rotor and to describe the action of
controls of the system. In addition, Wang in \cite{wa13e} gives the
regular point reduction by stages and Hamilton-Jacobi theorem of
regular controlled Hamiltonian system with symmetry on the
generalization of a semidirect product Lie group, and as an
application of the theoretical result, she studies the underwater
vehicle with two internal rotors as a regular point reducible by
stages RCH system, by using semidirect product Lie group and
Hamiltonian reduction by stages. On the other hand, if we define a
controlled Hamiltonian system on the cotangent bundle $T^*Q$ by
using a Poisson structure, see Wang and Zhang in \cite{wazh12} and
Ratiu and Wang in \cite{rawa12}, and the way of symplectic reduction
for regular controlled Hamiltonian system cannot be used, what and
how we could do? This is a problem worthy to be considered in
detail. In addition, we also note that there have been a lot of
beautiful results of reduction theory of Hamiltonian systems in
celestial mechanics, hydrodynamics and plasma physics. Thus, it is
an important topic to study the application of reduction theory of
controlled Hamiltonian systems in celestial mechanics, hydrodynamics
and plasma physics. These are our goals in future research.\\

\noindent{\bf Acknowledgments:} The author would like to thank
Professor Ivailo Mladenov for his interest and help for the research
work, and also to thank the support of Nankai University, 985
Project and the Key Laboratory of Pure Mathematics and
Combinatories, Ministry of Education, China.\\

\end{document}